\documentclass[12pt,a4paper]{article}

\usepackage{amsmath,amsfonts,amssymb} 
\allowdisplaybreaks
\newtheorem{theorem}{Theorem}
\newtheorem{Lem}{Lemma}

\newenvironment{pf}{\noindent{\bf Proof.}\enspace}{
}

\def\1{\mbox{1\hspace{-.25em}I}}
\def\Pb{\mathbf{P}}

\def\Ex{\mathbf{E}}

\DeclareMathOperator*{\argmin}{arg\,min}

\begin{document}

\title{On goodness-of-fit tests for parametric hypotheses in perturbed dynamical systems using a minimum distance estimator
}

\author{Maroua BEN ABDEDDAIEM\\
{\small Laboratoire Manceau de Math\'ematiques, Universit\'e du Maine}\\
 Le Mans,  France
}

\date{}


\maketitle

\begin{abstract}
We consider the problem of the construction of the Goodness-of-Fit test in the
case of continuous time observations of a diffusion process with
small noise.  The null hypothesis is parametric and we use a minimum distance
estimator of the unknown parameter. We propose an asymptotically
distribution free test for this model.
\end{abstract}

\section{Introduction}
\label{ch1I}
We consider the problem of the construction of the asymptotically distribution
free (ADF) test in the case of continuous time observations of a diffusion process. We suppose that
under the null hypothesis the trend coefficient depends on some unknown
one-dimensional parameter. Therefore the basic (null) hypothesis is
parametric.

The goodness-of-fit (GoF) test plays an important role in mathematical
statistics because its application allows to check if the observations
correspond well to the proposed mathematical model. Remind that in the case of
i.i.d. observations the problem of the construction of GoF tests attracts
attention of the statisticians since 1928 due to the works by Cram\'er (1928), von
Mises (1931) and Smirnov (1937) (see, e.g., \cite{Dar57} and \cite{Du}).

Let us recall the well-known basic results for the i.i.d. model. Suppose
that the null hypothesis is simple. Denoting the continuous distribution function under the null hypothesis by $F_0\left(x\right)$, we have to
check if the i.i.d. observations $X^n=\left(X_1,\ldots,X_n\right)$ have this
continuous distribution function. Many GoF tests are based on the following
property: the normalized empirical distribution function
$$\sqrt{n} \left( \hat F_n\left(x\right)-F_0\left(x\right)\right), \qquad
\hat{F}_n\left(x\right)=\frac{1}{n}\sum_{j=1}^n \1_{\{X_j < x\}}$$  converges in distribution, under the
null hypothesis, to the Brownian bridge
$B\left(F_0\left(x\right)\right)$.
In particular, for the
Cram\'er-von Mises statistic we have (with the change of variable
$s=F_0(x)$)
\begin{align*}
\Delta _n&={n}\int_{-\infty }^{\infty }
\left( \hat F_n\left(x\right)-F_0\left(x\right)\right)^2{\rm
  d}F_0\left(x\right) \\
&\qquad \Longrightarrow  \int_{-\infty }^{\infty
}B\left(F_0\left(x\right)\right)^2{\rm d} F_0\left(x\right)=\int_{0}^{1}B\left(s\right)^2{\rm d}s\equiv \Delta .
\end{align*}
Therefore the limit $\Delta $ does not depend on $F_0\left(\cdot \right)$
(distribution free) and
the Cram\'er-von Mises test
\begin{equation*}
\psi_n=\1_{\left\{ \Delta _n>c_\alpha \right\}},\qquad \Pb\left\{\Delta >c_\alpha \right\}=\alpha
\end{equation*}
has the asymptotic $(n\rightarrow \infty)$ size $\alpha \in \left(0,1\right)$ (see, e.g., \cite{LR}).

 In the case of
parametric null hypothesis
$$\mathcal{H}_0 \quad: \qquad F\left(x\right)=F_0\left(\vartheta ,x\right),\qquad \vartheta\in \Theta \subseteq \mathbb{R},$$
where $\vartheta$
 is some unknown one-dimensional parameter, the similar limit
\begin{align*}
U_n\left(x\right) = \sqrt{n}\left(\hat F_n\left(x\right)-F_0(\hat\vartheta
_n,x)\right)
\Longrightarrow  B\left(F_0\left(\vartheta,x\right)\right)-\zeta \dot
F_0\left(\vartheta ,x\right)\equiv U\left(x\right)
\end{align*}
 is no more
distribution free and the choice of the threshold $c_\alpha $ is much more
complicated. Here $\hat\vartheta _n$ is the maximum likelihood estimator (MLE), $\dot{F}_0\left(\vartheta
  ,\cdot\right)$ is the derivative of $F_0\left(\vartheta ,\cdot\right)$
  w.r.t. $\vartheta$ and $\zeta$ is a Gaussian variable.
We have for $U\left(x\right)$, by the change of variables $t=F_0\left(\vartheta ,x\right)$ and $s=F_0\left(\vartheta ,y\right)$,
  the representation (see, e.g., Darling
  \cite{Dar57})
\begin{equation}
\label{ch101}
u\left(t\right)=B\left(t\right)-\int_{0}^{1}h\left(\vartheta,s\right){\rm
  d}B\left(s\right)\int_{0}^{t}h\left(\vartheta,s\right){\rm d}s,
\end{equation}
where
$$
h\left(\vartheta,s\right)=\frac{\dot{l}\left(\vartheta,F_\vartheta^{-1}(s)\right)}
{\sqrt{I\left(\vartheta\right)}},\qquad\int_{0}^{1}h\left(\vartheta,s\right)^2{\rm d}s=1.
$$
Here $l\left(\vartheta,x\right)=\ln f\left(\vartheta,x\right),$ $\dot{l}\left(\vartheta,x\right)$ is the derivative of $l\left(\vartheta,x\right)$ w.r.t. $\vartheta$,
$f\left(\vartheta,x\right)$ is the density function, $I\left(\vartheta\right)$
is the Fisher information and $y=F_\vartheta^{-1}(s)$ is the inverse function to
$F_0\left(\vartheta ,y\right)$, i.e., $F_0(\vartheta ,y)=s$. Note that the function $h\left(\cdot,\cdot
\right)$ appears here twice because we used the MLE.

 For the minimum distance estimator (MDE) using the same arguments as in the
   case of the MLE we obtain the limit
\begin{equation}
\label{ch102}
u\left(t\right)=B\left(t\right)-\int_{0}^{1}g\left(\vartheta,s\right)\;{\rm
  d}B\left(s\right)\int_{0}^{t}h\left(\vartheta,s\right)\;{\rm d}s,\qquad
\int_{0}^{1}g\left(\vartheta,s\right)^2{\rm d}s=1
\end{equation}
with two different functions. In both cases, the corresponding tests are no more
asymptotically distribution free and the choice of the
threshold from the equation
\begin{align*}
\Pb\left\{ \int_{0}^{1}u\left(s\right)^2{\rm d}s>c_\alpha \right\}=\alpha
\end{align*}
can be a difficult problem.

There are several possibilities to solve this problem in the case of the limit
\eqref{ch101}. One of them is to find a linear transformation $L\left[\cdot
  \right]$ of the random function $u\left(\cdot \right)$, such that
$L\left[u\right](t)=W\left(t\right)$, where
$W\left(t \right),0\leq t\leq 1$ is some Wiener process. Then the statistic
\begin{align*}
\hat\Delta _n&=\int_{-\infty }^\infty L\left[ U_n\right]\left(x\right)^2{\rm
  d}F_0(\hat\vartheta _n,x) \Longrightarrow
\int_{-\infty }^\infty W\left(F_0\left(\vartheta ,x\right)\right)^2{\rm d} F_0\left(\vartheta
,x\right)\\
&\qquad =\int_{0}^{1}W\left(t\right)^2{\rm d}t \equiv \hat{\Delta}
\end{align*}
and we obtain once more the distribution free limit $\hat\Delta$. Therefore the
test
$$
\hat\psi _n=\1_{\left\{\hat\Delta _n>d_\alpha \right\}},\qquad\qquad
\Pb\left\{\hat{\Delta} >d_\alpha \right\}=\alpha
$$
is ADF. Such transformation was proposed by Khmaladze \cite{Kh81}. We have to emphasize that in \cite{Kh81} and in many other works (see, e.g., the paper Maglapheridze et al. \cite{MTP}) the estimator used was always the MLE and this is important for the construction of this linear transformation. Many authors wrote that similar transformation can be obtained in the case of other estimators with limit \eqref{ch102}, but as we know this work (construction of the linear transformation with other estimators) was not done.

The problem of GoF testing for the model of continuous
time observations of diffusion process, with a simple null hypothesis
$\Theta=\{\vartheta_{0}\}$, was studied in \cite{DK} and \cite{K11}.
 Suppose that the observed diffusion process under hypothesis is
\begin{equation}
\label{ch1}
\textrm{d}X_t=S_0(X_t)\;\textrm{d}t+\varepsilon\sigma\left(X_t\right)\;\textrm{d}W_t,\qquad X_0=x_0, \quad  0\leq t \leq T,
\end{equation}
with deterministic initial value $x_0$, known diffusion coefficient
$\varepsilon^2 \sigma\left(\cdot\right)^2>0$, some known smooth function
$S_0(\cdot)$ and $W_t, 0 \leq t \leq 1$ is a Wiener process.
Note that if the functions $S_0\left(x\right)$ and $\sigma\left(x\right)$ are Lipschitz w.r.t. $x$, then we have with probability 1
$$
\sup_{0\leq t \leq T}\left|X_t-x_t \right| \leq C\varepsilon \sup_{0\leq t \leq T}\left|W_t\right|,
$$
where $C$ is some constant. Therefore the process $X_t$ converges to $x_t=x_t(\vartheta_0)$ (solution of the
equation \eqref{ch1} as $\varepsilon =0$) uniformly w.r.t. $t\in\left[0,T\right]$ with probability 1. We have as well
$$
\sup_{0\leq t \leq T}\Ex \left|X_t-x_t \right|^2 \leq C\varepsilon^2.
$$
For the proof see  \cite{Kut94}, Lemma 1.13.

 The GoF test was constructed on the basis of the normalized difference $\varepsilon^{-1}\left(X_t-x_t\right)$ and the limit of this statistic is a Gaussian process. This process can be transformed into a Wiener process $w(s), 0 \leq s \leq 1$ as follows.

Introduce the statistic
$$
\delta_\varepsilon = \left[\int_0^T\left(\frac{\sigma\left(x_t\right)}{S_0(x_t)}\right)^2\textrm{d}t\right]^{-2}\int_0^T
\left(\frac{X_t-x_t}{\varepsilon S_0(x_t)^2}\right)^2\sigma\left(x_t\right)^2\textrm{d}t.
$$
The following convergence
$$
\delta_\varepsilon \Longrightarrow \delta \equiv \int_0^1w(s)^2\;\textrm{d}s
$$
was proved and therefore the test $\phi_\varepsilon = \1_{\{\delta_\varepsilon > c_\varepsilon\}}$ with
$\Pb\left(\delta>c_\varepsilon\right)=\varepsilon$ is ADF.

The case
of parametric basic hypothesis and ADF tests for ``small noise" diffusion
processes was studied, for example, in \cite{Ikut01}, \cite{Kut13b}$-$\cite{Kut13c}
and the estimator used, for the construction of the linear
transformation and tests, was always the MLE (see Kutoyants \cite{Kut78}). There are several ADF GoF tests
for the ergodic diffusion processes proposed, for example, in the works
\cite{NN}$-$\cite{NZ} and \cite{KK14}.

Note that in some problems it is not possible to have an explicit
  expression for the MLE and therefore  sometimes it is better to use other
  estimators, which can be easily calculated. For example, this can be the minimum
  distance estimator, a method of
  moments estimator or a trajectory fitting estimator and so on. In all such cases, the limit
  expression for the underlying statistics will be like \eqref{ch102} but with
  two different functions $h\left(\cdot,\cdot\right)$ and
  $g\left(\cdot,\cdot\right)$ (see below).

In this work, we observe continuous time process $X^{\varepsilon}$, which is the solution of the stochastic differential equation
$$
\textrm{d}X_t=S(X_t)\;\textrm{d}t+\varepsilon\;\textrm{d}W_t,\qquad X_0=x_0, \quad  0\leq t \leq T,
$$
where $W_t, \;0\leq t\leq T$ is a Wiener process and $S(x)$ is some unknown smooth function. Based on a minimum distance estimator, we construct an ADF test in the case of parametric null hypothesis
$$
\mathcal{H}_0: \quad S(x) = S(\vartheta,x),\qquad  \vartheta \in
\Theta=(a,b).
$$

The main results  are presented in Theorem \ref{ch1T1} and Theorem
  \ref{ch1T2}. In particular, the ADF test is given in Theorem \ref{ch1T2} below.
We realize the following program. First, we show that the basic statistic
$$
u_\varepsilon \left(t\right)=\frac{X_t-x_t\left(\vartheta _\varepsilon
  ^*\right)}{\varepsilon \;S\left(\vartheta _\varepsilon
  ^*,X_t \right)}, \qquad \qquad 0\leq t \leq T
$$
($\vartheta _\varepsilon^*$ is the MDE) converges to the random process $u\left(t\right), 0\leq t \leq T$ (see \eqref{ch1001} below). Then we transform $u\left(t\right)$ in $U\left(\frac{t}{T} \right), 0\leq t \leq T$ (see \eqref{ch1002} below). We obtain for $U\left(\cdot \right)$ the following representation
\begin{equation}
\label{ch103}
U\left(\nu\right)=W\left(\nu\right)-\int_{0}^{1}g\left(\vartheta,r\right){\rm
  d}W\left(r\right)\int_{0}^{\nu}h\left(\vartheta,r\right){\rm d}r,\quad
\int_{0}^{1}g\left(\vartheta,r\right)^2{\rm d}r=1.
\end{equation}
The last step is to apply the special linear transformation $L\left[U\right]\left(\nu \right)=w_\nu , 0\leq \nu \leq 1$, where $w_\nu, 0\leq \nu \leq 1$ is a Wiener process (see Theorem
\ref{ch1T1}). This allows us to construct an ADF test (see Theorem \ref{ch1T2}). The main contribution of this work is the form of this linear transformation.

Then  we realize similar transformations with the ``empirical''
process $u_\varepsilon \left(\cdot\right)$, apply the linear transformation
$L\left[\cdot \right]$ and obtain the convergence
$$
\Delta _\varepsilon =\frac{1}{T}\int_{0}^{T}L\left[U_\varepsilon \right]\left(t
\right)^2{\rm d}t \Longrightarrow \int_{0}^{1}w_\nu ^2\;{\rm d}\nu.
$$
Here the process $U_\varepsilon(\cdot)$ will be defined in Section 5 by \eqref{ch128}.
Therefore the test $\psi _\varepsilon =\1_{\left\{\Delta _\varepsilon>c_\alpha
  \right\}}$
 will be ADF, because the limit distribution of $\Delta _\varepsilon $ does
not depend on $S\left(\cdot \right)$ and $\vartheta$ and the test is of asymptotic
$(\varepsilon \rightarrow 0)$ size $\alpha\in(0,1)$.

Note that if in our problem we use the MLE of the unknown parameter, then the
limit representation is
\begin{equation}
\label{ch104}
U\left(\nu\right)=W\left(\nu\right)-\int_{0}^{1}h\left(\vartheta,s\right){\rm
  d}W\left(s\right)\int_{0}^{\nu}h\left(\vartheta,s\right){\rm d}s,\qquad
\int_{0}^{1}h\left(\vartheta,s\right)^2{\rm d}s=1
\end{equation}
and  our transformation coincides with the one proposed in \cite{Kh81}. We
discuss this case in Section 6.

We have to note that this linear transformation is rather cumbersome and the
realization of the test can be a computationally difficult problem too. We
suppose that the presented result is of theoretical interest and allows  to ``close the
gap'' in this field.

At the same time, we understand that this result is in some sense ``negative'' and says that if we have no MLE it is better to seek another GoF test, which is ADF. Note as well that in i.i.d. case, even if the estimated parameter is one-dimensional, the reduction of the equation with Brownian bridge \eqref{ch102} to the equation \eqref{ch103}, using the relation $(B\left(t\right)=W(t) -t \;W(1))$, leads to the corresponding Fredholm equation (see \eqref{ch116} below) with two-dimensional $g\left(\cdot\right)$ and $h\left(\cdot\right)$. The expression for the solution of this equation and the form of the linear transformation becomes much more complicated. This is probably the reason why this problem was not considered till now. Our results thus should be understood as a constructive existence result for ADF tests based on the MDE.

\section{Minimum distance estimator}
\label{ch1s01}
 Suppose that the continuous time observed process
$X^{\varepsilon} = (X_t,0\leq t \leq T)$ is the solution of the stochastic
differential equation
\begin{equation}
\label{ch105}
\textrm{d}X_t=S(X_t)\;\textrm{d}t+\varepsilon\;\textrm{d}W_t,\qquad X_0=x_0, \quad  0\leq t \leq T,
\end{equation}
where $W_t, \;0\leq t\leq T$ is a Wiener process, the initial value $x_0$ is
deterministic and the trend coefficient $S(x)$ is some unknown smooth function.

We consider the composite basic hypothesis
$$
\mathcal{H}_0: \quad S(x) = S(\vartheta,x),\qquad  \vartheta \in
\Theta=(a,b),
$$
 where $\vartheta$ is the one-dimensional (unknown) parameter, against
alternative $\mathcal{H}_1:$ not $\mathcal{H}_0$, i.e., the trend
coefficient $S\left(x\right)$ in the observed diffusion process \eqref{ch105}
does not belong to the parametric family $\left\{S\left(\vartheta
,x\right),\vartheta \in \Theta \right\}$. Therefore the process, under
hypothesis $\mathcal{H}_0$, has the stochastic differential
\begin{equation}\label{ch1s00}
\textrm{d}X_t=S(\vartheta,X_t)\;\textrm{d}t+\varepsilon\;\textrm{d}W_t,\qquad
X_0=x_0, \quad 0\leq t \leq T.
\end{equation}

We are interested in the properties of the test in the asymptotics of {\it
small noise}, i.e., as $\varepsilon\rightarrow 0$.
Below and in the sequel the dot means derivation w.r.t. $\vartheta$.

Let us introduce the regularity conditions.\\ ${\cal R}.$ {\it The function
  $S(\vartheta,x)$ is strictly positive and has two continuous bounded
  derivatives with respect to $\vartheta$ and $x$.}

In the presentation below we suppose that these conditions and the
  basic hypothesis $\mathcal{H}_0$ are always fulfilled.

It is known that the
solution $X_t$ converges uniformly in $t\in [0,T]$ to the solution $x_t =
x_t(\vartheta)$ of the ordinary differential equation
$$\frac{\textrm{d}x_t}{\textrm{d}t} = S(\vartheta,x_t),\qquad x_0 ,\quad 0\leq t \leq T,$$
where $x_0$ is the same as in \eqref{ch105} (see \cite{FW} or \cite{Kut94}, Lemma
1.13 for the proof).

Recall the properties of  maximum likelihood and minimum distance estimators
of the parameter $\vartheta $. The
likelihood ratio function in the case of observations \eqref{ch1s00} is
$$
\ell\left(\vartheta ,X^\varepsilon\right)=\frac{\textrm{d}P_\vartheta}{\textrm{d}P_0}=\exp\left\{\int_0^T\frac{S\left(\vartheta,X_t\right)}{\varepsilon^2}
\;\textrm{d}X_t-\int_0^T\frac{S\left(\vartheta,X_t\right)^2}{2\;\varepsilon^2}\;\textrm{d}t\right\}, \quad \vartheta \in \Theta,
$$
where $P_\vartheta$ and $P_0$ are the measures induced respectively by the processes \eqref{ch1s00} and
$$
\textrm{d}X_t=\varepsilon \;\textrm{d}W_t, \qquad X_0=x_0, \qquad 0 \leq t \leq T,
$$
(see \cite{LS} for more details).

The MLE $\hat\vartheta_\varepsilon $ is solution of the
equation
\begin{align*}
\sup_{\vartheta \in \Theta }\ell\left(\vartheta ,X^\varepsilon
\right)=\ell(\hat\vartheta_\varepsilon  ,X^\varepsilon ) .
\end{align*}
This estimator is consistent and asymptotically normal (as $\varepsilon
\rightarrow 0$)
\begin{align*}
\varepsilon ^{-1}(\hat\vartheta_\varepsilon-\vartheta)\Longrightarrow {\cal N}\left(0,{\rm I}\left(\vartheta
\right)^{-1}\right),\qquad {\rm I}\left(\vartheta \right)=\int_{0}^{T}\dot
S\left(\vartheta ,x_t\right)^2{\rm d}t.
\end{align*}
Here and in the sequel $x_t=x_t\left(\vartheta\right)$. For the proof see, Kutoyants \cite{Kut94}.

Introduce the {\it minimum distance estimator}
$$\vartheta^{*}_{\varepsilon} =\argmin\limits_{\vartheta\in \Theta} \left\|X-x(\vartheta)\right\| ,$$
where $\|\cdot\|$ is $L^2[0,T]$ norm defined by
$$
\left\| X-x(\vartheta)\right\|^2 = \int_0^T
\left(X_t-x_t(\vartheta)\right)^2\textrm{d}t.
$$
The properties of the MDE for this model were studied in \cite{Kut94}.
The MDE satisfies the {\it minimum distance equation} (MDEq)
$$
\int_0^T \left(X_t-x_t(\vartheta^{*}_{\varepsilon})\right)\dot{x}_t(\vartheta^{*}_{\varepsilon})\;\textrm{d}t=0.
$$
Let us put $u_{\varepsilon}^{*}=\varepsilon
^{-1}\left(\vartheta^{*}_{\varepsilon}-\vartheta \right)$. Therefore the MDEq is as follows:
$$
\int_0^T \left(X_t-x_t\left(\vartheta + \varepsilon
u_{\varepsilon}^{*}\right)\right)\dot{x}_t(\vartheta^{*}_{\varepsilon})\;\textrm{d}t=0.
$$
Then,  by the Taylor
formula, we can write
$$
\int_0^T \left(X_t- x_t(\vartheta) - \varepsilon
\;u_{\varepsilon}^{*} \;\dot{x}_t(\tilde{\vartheta})\right)\dot{x}_t(\vartheta^{*}_{\varepsilon})\;\textrm{d}t=0,
$$
where
$ |\tilde\vartheta -\vartheta |\leq \left| \vartheta^{*}_{\varepsilon}-\vartheta \right|$
and for $u_\varepsilon ^*$ we obtain the following representation
\begin{equation}
\label{ch11}
u_{\varepsilon}^{*}  =
J(\vartheta)^{-1}\int_0^T\varepsilon^{-1}\left(X_t-x_t(\vartheta)\right)\dot{x}_t(\vartheta)\;\textrm{d}t+o(1),
\end{equation}
where
\begin{equation}
\label{ch12}
J(\vartheta) =
\int_0^T\dot{x}_t(\vartheta)^2\textrm{d}t.
\end{equation}
Let us consider the random process $x_t^{(1)}$ defined as the derivative of $X_t$ w.r.t.
$\varepsilon$ at $\varepsilon=0$, i.e., we have $\varepsilon
^{-1}\left(X_t-x_t\left(\vartheta \right)\right)\rightarrow x_t^{(1)}$. Here the
random process $x_t^{(1)}$ satisfies the linear equation
$$
\textrm{d}x_t^{(1)} = S'(\vartheta,x_t)\;x_t^{(1)}\textrm{d}t +
\textrm{d}W_t,\qquad x_0^{(1)}=0,\quad0\leq t \leq T,
$$
where
$S'(\vartheta,x) = \displaystyle\frac{\partial S(\vartheta,x)}{\partial x}.$
The process $x_t^{(1)}$ has the representation
\begin{equation}
\label{ch1011}
x_t^{(1)}=S(\vartheta,x_t)\int_0^t \frac{1}{S(\vartheta,x_s)}\;\textrm{d}W_s
\end{equation}
(see, e.g., Section $3.3$ in \cite{Kut94}).
Due to the above representation of $x_t^{(1)},$ the equation \eqref{ch11} becomes
\begin{eqnarray*}
u_{\varepsilon}^{*}&
= & J(\vartheta)^{-1}\int_0^Tx_t^{(1)}\dot{x}_t(\vartheta)\;\textrm{d}t+o(1)\\ & = &
J(\vartheta)^{-1}\int_0^TS(\vartheta,x_t)\int_0^t
\frac{1}{S(\vartheta,x_s)}\;\textrm{d}W_s\;\dot{x}_t(\vartheta)\;\textrm{d}t+o(1).
\end{eqnarray*}
Here $o(1)$ means the convergence in probability, i.e., for any $\nu>0$, we have
\begin{align}
\label{cch10}
\lim_{\varepsilon \rightarrow 0}\Pb_\vartheta \left(\left|u_{\varepsilon}^{*}-\rho(\vartheta)\right|>\nu\right)=0,
\end{align}
where
$$
\rho(\vartheta)=J(\vartheta)^{-1}\int_0^TS(\vartheta,x_t)\int_0^t
\frac{1}{S(\vartheta,x_s)}\;\textrm{d}W_s\;\dot{x}_t(\vartheta)\;\textrm{d}t.
$$

Therefore $\vartheta^{*}_{\varepsilon}$ admits the following representation, by Fubini's Theorem,
\begin{equation}
\label{ch13}
\varepsilon^{-1}\left(\vartheta^{*}_\varepsilon-\vartheta\right)=J(\vartheta)^{-1}\int_0^T\frac{1}{S(\vartheta,x_v)}\int_v^T S(\vartheta,x_s)\;\dot{x}_s(\vartheta)\;\textrm{d}s\;\textrm{d}W_v+o(1).
\end{equation}

Moreover, under conditions of regularity, the estimator $\vartheta^{*}_{\varepsilon}$
is consistent and asymptotically normal (see Chapter 7 in \cite{Kut94})
$$\mathcal{L}_{\vartheta}\{\varepsilon^{-1}(\vartheta^{*}_{\varepsilon}-\vartheta)\}\Longrightarrow \mathcal{L}\{\xi\}=\mathcal{N}(0,\sigma^2\left(\vartheta \right)),$$
where
$$
\sigma^2\left(\vartheta \right) = J(\vartheta)^{-2}\int_0^T\frac{1}{S(\vartheta,x_v)^2}\left(\int_v^TS(\vartheta,x_s)\;\dot{x}_s(\vartheta)\;\textrm{d}s\right)^2\textrm{d}v.
$$

\section{Basic statistic}
\label{ch1s02}

 Our goal is to find the GoF test, which is ADF, i.e., we seek the test
 statistic whose limit distribution, under hypothesis, does not depend on the
 underlying model given by the function $S(\vartheta,x)$ and parameter
 $\vartheta$.

  Introduce the statistic
$$
\delta^*_{\varepsilon}=\int_0^T
\left[\frac{X_t-x_t(\vartheta^{*}_{\varepsilon})}{\varepsilon \;
    S\left(\vartheta^{*}_{\varepsilon},X_t\right)}\right]^2\textrm{d}t.
$$
To study it we need the behavior of the difference
$X_t-x_t(\vartheta^{*}_{\varepsilon})$, which can be described as follows:
\begin{eqnarray*}
&&\varepsilon^{-1}\left(X_t-x_t(\vartheta^{*}_{\varepsilon})\right) =
  \varepsilon^{-1}\left(X_t-x_t(\vartheta)\right)-\varepsilon^{-1}
  \left(\vartheta^{*}_{\varepsilon}-\vartheta\right)\dot{x}_t(\vartheta)+o(1)\\ &&\qquad
  \qquad =
  x_t^{(1)}(\vartheta)-J(\vartheta)^{-1}\int_0^Tx_s^{(1)}\;\dot{x}_s
(\vartheta)\;\textrm{d}s\;\dot{x}_t(\vartheta)+o(1)\\ &&\qquad
  \qquad = S(\vartheta,x_t)\int_0^t
  \frac{\textrm{d}W_v}{S(\vartheta,x_v)}\\ && \qquad\quad \quad
  -J(\vartheta)^{-1}\int_0^T\frac{1}{S(\vartheta,x_v)}\int_v^T
  S(\vartheta,x_s)\;\dot{x}_s(\vartheta)\;\textrm{d}s\;\textrm{d}W_v\;\dot{x}_t(\vartheta)+o(1),
\end{eqnarray*}
where the process $x_t^{(1)}$ was defined by \eqref{ch1011} and the derivative
$\dot{x}_t(\vartheta)$ w.r.t. $\vartheta$ satisfies the equation
$$
\frac{\textrm{d}\dot{x}_t(\vartheta)}{\textrm{d}t} =
S'(\vartheta,x_t) \;\dot{x}_t(\vartheta)+\dot{S}(\vartheta,x_t),
\qquad\dot{x}_0(\vartheta)=0.
$$
Its solution is the function (it can be found in \cite{Kut13b})
$$
\dot{x}_t(\vartheta) = S(\vartheta,x_t)\int_0^t
\frac{\dot{S}(\vartheta,x_v)}{S(\vartheta,x_v)}\;\textrm{d}v.
$$
Here $o(1)$ is the uniform convergence in probability w.r.t. $t\in[0,T]$, i.e., for any $\nu>0$, we have
\begin{align}
\label{cch11}
\lim_{\varepsilon \rightarrow 0}\Pb_\vartheta \left(\sup_{t\in[0,T]} \left|\varepsilon^{-1}\left(X_t-x_t(\vartheta^{*}_{\varepsilon})\right)-\mu_t(\vartheta)\right|>\nu\right)=0,
\end{align}
where
\begin{align*}
&\mu_t(\vartheta)=S(\vartheta,x_t)\int_0^t
  \frac{\textrm{d}W_v}{S(\vartheta,x_v)}
  \\&\qquad\qquad-J(\vartheta)^{-1}\int_0^T\frac{1}{S(\vartheta,x_v)}\int_v^T
  S(\vartheta,x_s)\;\dot{x}_s(\vartheta)\;\textrm{d}s\;\textrm{d}W_v\;\dot{x}_t(\vartheta).
\end{align*}
For the details see \cite{Kut94}, Chapter 7.

Hence we have the uniform convergence w.r.t. $t\in[0,T]$ (in probability)
\begin{align}
\begin{split}
\label{ch1001}
&u_{\varepsilon}(t) =
  \frac{X_t-x_t(\vartheta^{*}_{\varepsilon})}{\varepsilon\;S\left(\vartheta^{*}_{\varepsilon},X_t\right)}
  \longrightarrow u(t) = \int_0^t
  \frac{\textrm{d}W_v}{S(\vartheta,x_v)}\\ &\qquad-J(\vartheta)^{-1}\int_0^T\frac{1}{S(\vartheta,x_v)}
  \int_v^TS(\vartheta,x_s)\;\dot{x}_s(\vartheta)\;\textrm{d}s\;\textrm{d}W_v\int_0^t
  \frac{\dot{S}(\vartheta,x_v)}{S(\vartheta,x_v)}\;\textrm{d}v
  \end{split}
\end{align}
and it can be shown (see  in \cite{Kut13b} the details in the similar problem where the MLE was used)
 that
$$\delta^*_{\varepsilon}\Longrightarrow \int_0^Tu(t)^2\;\textrm{d}t.$$
Therefore the test based on this statistic is not ADF. Hence to obtain an ADF GoF
test we introduce the Gaussian process
\begin{equation}
\label{ch1002}
U\left(\frac{t}{T}\right)  =  \frac{1}{\sqrt{T}}\int_0^{t}
S(\vartheta,x_s)\;\textrm{d}u(s),\qquad 0\leq t\leq T,
\end{equation}
where $u\left(\cdot\right)$ was given by \eqref{ch1001}. Then, by It\^{o} formula,
\begin{align*}
\textrm{d}\left(S(\vartheta,x_s)u(s)\right)=\textrm{d}S(\vartheta,x_s)\;u(s)
+\textrm{d}u(s)\;S(\vartheta,x_s)
\end{align*}
and using the equality
\begin{equation}
\label{ch1005}
\textrm{d}S(\vartheta,x_s)=S'(\vartheta,x_s)\;S(\vartheta,x_s)\;\textrm{d}s,
\end{equation}
we have
$$
S(\vartheta,x_t)\;u(t)=\int_0^tS'(\vartheta,x_s)\;S(\vartheta,x_s)\;u(s)\;\textrm{d}s
+\int_0^tS(\vartheta,x_s)\;\textrm{d}u(s).
$$
Here and in the sequel we denoted by prime the derivative w.r.t. $x$.

Therefore the process $U(\cdot)$ defined by \eqref{ch1002} admits the following representation $(0\leq t \leq T)$
\begin{align}
\label{ch10004}
&U\left(\frac{t}{T}\right)  =  \frac{1}{\sqrt{T}}\;S(\vartheta,x_t)\;u(t)
- \frac{1}{\sqrt{T}}
\int_0^tS'(\vartheta,x_s)\;S(\vartheta,x_s)\;u(s)\;\textrm{d}s.
\end{align}

Further we define two functions
\begin{align}
\begin{split}
\label{ch107}
h(\vartheta,r)
&=T\tilde{J}(\vartheta)^{-1}\dot{S}(\vartheta,x_{rT}(\vartheta))\;C\left(\vartheta\right)^{1/2}
\end{split}
\end{align}
and
\begin{align}
\begin{split}
\label{ch108}
g(\vartheta,r) &= S(\vartheta,x_{rT}
(\vartheta))^{-1}\displaystyle\int_{r}^1S(\vartheta,x_{zT}(\vartheta))\;\dot{x}_{zT}(\vartheta)\;\textrm{d}z
\;C\left(\vartheta\right)^{-1/2}.
\end{split}
\end{align}
Here
$
\tilde{J}(\vartheta) = \displaystyle\int_0^1 \dot{x}_{vT}(\vartheta)^2 \textrm{d}v$ and
$$
C\left(\vartheta\right)=\int_0^1S(\vartheta,x_{vT}(\vartheta))^{-2}\left(\int_{v}^1S(\vartheta,x_{zT}(\vartheta))
\;\dot{x}_{zT}(\vartheta)\;\textrm{d}z\right)^2\textrm{d}v.
$$
Observe that
\begin{align*}
\int_0^1 g(\vartheta,r)^2\textrm{d}r=1.
\end{align*}

\begin{Lem} \label{ch1L1} We have the equality
\begin{equation}
\label{ch109}
U(\nu)  =  W(\nu)  - \int_0^1g(\vartheta,r)\;\textrm{d}W(r)\;\int_0^\nu h(\vartheta,r)\;\textrm{d}r, \quad 0\leq \nu \leq 1,
\end{equation}
where $W(\nu), 0\leq\nu\leq1$ is a Wiener process.
\end{Lem}

\begin{pf}
We have by It\^{o} formula, for $0\leq t \leq T$,
\begin{align*}
&\textrm{d}U\left(\frac{t}{T}\right)=
\frac{1}{\sqrt{T}}\;\textrm{d}S(\vartheta,x_t)\;u(t)\\
&\qquad \qquad\quad+\frac{1}{\sqrt{T}}\;S(\vartheta,x_t)\;\textrm{d}u(t)- \frac{1}{\sqrt{T}}\;
S'(\vartheta,x_t)\;S(\vartheta,x_t)\;u(t)\;\textrm{d}t.
\end{align*}
Then by \eqref{ch1005} and using the representation \eqref{ch1001} of the process
$u(\cdot)$, we obtain
\begin{align*}
&U\left(\frac{t}{T}\right)=\frac{1}{\sqrt{T}}\int_0^tS(\vartheta,x_s)\;\textrm{d}u(s)\\
&=\frac{W_t}{\sqrt{T}}-
\frac{J(\vartheta)^{-1}}{\sqrt{T}} \int_0^T\frac{1}{S(\vartheta,x_v)}
\int_v^TS(\vartheta,x_s)\;\dot{x}_{s}(\vartheta)\;\textrm{d}s\;\textrm{d}W_v
\int_0^t\dot{S}(\vartheta,x_v)\;\textrm{d}v.
\end{align*}
Hence, by the change of variables $\nu=\displaystyle\frac{t}{T}$ and $W(\nu) = T^{-1/2}W_{\nu T}, 0\leq \nu \leq 1$, we have
\begin{align*}
&U(\nu)  =
W\left(\nu\right) \\
&\qquad\qquad-
\frac{J(\vartheta)^{-1}}{\sqrt{T}} \int_0^T\frac{1}{S(\vartheta,x_v)}
\int_v^TS(\vartheta,x_s)\;\dot{x}_{s}(\vartheta)\;\textrm{d}s\;\textrm{d}W_v\int_0^{\nu
T}\dot{S}(\vartheta,x_v)\;\textrm{d}v.
\end{align*}
Let us change the variables
$r=\displaystyle\frac{v}{T}$,
$z=\displaystyle\frac{s}{T}$, $W(r) = T^{-1/2}W_{rT}, 0\leq r \leq 1$. Then we can write
\begin{align*}
U(\nu)& = W(\nu) \\
& \quad- \tilde{J}(\vartheta)^{-1} \int_0^1
\frac{1}{S(\vartheta,x_{rT})} \int_{rT}^TS(\vartheta,x_s)\;\dot{x}_s(\vartheta)\;\textrm{d}s\;\textrm{d}W(r)\int_0^{\nu
}\dot{S}(\vartheta,x_{rT})\;\textrm{d}r\\
& = W(\nu) \\
& \quad-\frac{T}{\tilde{J}(\vartheta)}\int_0^1\frac{1}{S(\vartheta,x_{rT})} \int_{r}^1S(\vartheta,x_{zT})\;\dot{x}_{zT}(\vartheta)\;\textrm{d}z\;\textrm{d}W(r)\int_0^{\nu
}\dot{S}(\vartheta,x_{rT})\;\textrm{d}r\\
 & = W(\nu) -
\int_0^1g(\vartheta,r)\;\textrm{d}W(r)\;\int_0^{\nu}h(\vartheta,r)\;\textrm{d}r.
\end{align*}
Therefore we obtain the representation \eqref{ch109} and this proves the Lemma \ref{ch1L1}.
\end{pf}

It can be shown using the convergence of the \textit{empirical
    version} $U_\varepsilon\left(\cdot\right)$ to
  $U\left(\cdot\right)$ (see proof of Theorem \ref{ch1T2} below) and due to the
  continuous mapping Theorem, that we have the convergence
$$
\widetilde{\Delta}_{\varepsilon} = \frac{1}{T^2}\int_0^T\left( \int_0^{t}
S(\vartheta^*_{\varepsilon},X_s)\;\textrm{d}u_{\varepsilon}(s)\right)^2\textrm{d}t
\Longrightarrow \int_0^1 U(\nu)^2\textrm{d}\nu.
$$
We remark that the test based on this
statistic is not ADF. Hence we have to find the transformation $L[U](\cdot)$
into the Wiener process such that
$$
\int_0^1 L[U](\nu)^2\textrm{d}\nu= \int_0^1
w_{\nu}^2\;\textrm{d}\nu,\qquad 0\leq \nu\leq 1,
$$
where $w_{\nu}, 0\leq \nu\leq 1$ is a Wiener process.
This property allows us to present the ``empirical version'' of the test statistic
with the same limit. Therefore the test based on this statistic is ADF.

\section{Linear transformation}
\label{ch1s03}

Now the problem is to find such transformation $L\left[\cdot\right]$ of $U(\cdot)$ (see \eqref{ch109}) using the
  MDE  that $L[U](\nu)=w_{\nu}.$ Recall that for the limit
process \eqref{ch101} such linear transformation and the corresponding ADF test
were proposed by Khmaladze \cite{Kh81}. Another (direct) proof of this result
was recently obtained by Kleptsyna and Kutoyants \cite{KK14}. Note
  that in these works the estimator used was always the MLE and in our work it is the MDE. The limit processes in these two cases are quite different. That is why we have to present here a special modification of the proof given in \cite{KK14}.  Our proof follows the main
steps of the work \cite{KK14}. Specifically, we have to solve Fredholm
equation of the second kind with degenerated kernel. The solution of it gives
us the desired linear transformation.

Denote
\begin{align}
\begin{split}
&I_1 = \int_0^r g(\vartheta,q)^2\textrm{d}q,\; I_2 = \int_0^r h(\vartheta,q)
  \;g(\vartheta,q)\;\textrm{d}q,\; I_3 =\int_0^r
  h(\vartheta,q)\;\textrm{d}q,\\ &\qquad \qquad \quad I_4 =\int_0^r
  h(\vartheta,q)^2\textrm{d}q,\; I_5 =\int_0^r
  g(\vartheta,q)\;\textrm{d}q. \label{ch110}
  \end{split}
\end{align}
Below we omit $\vartheta$ and $r$ for simplicity and put $g=g(\vartheta,r)$
and $h=h(\vartheta,r)$. Introduce the functions
\begin{align}
\varphi_1(r) & =  g-h-3I_2g+I_5hg+  I_3g^2+2I_2h-2I_2I_3g^2+I_1I_2^2h+I_4I_5g^2-I_2^3g   \nonumber\\
& \quad -  I_2I_4g+3I_2^2g+  I_2I_5h^2-2I_2I_5hg-2I_1I_2h +  I_2^2I_5hg+I_1^2I_3h^2-I_4h\nonumber\\
&  \quad +  2I_1I_4h-I_1I_4g + I_1I_2I_4g + I_1I_4I_5hg-I_1^2I_4h+I_1h+2I_2I_3hg\nonumber\\
&  \quad -  I_2I_4I_5g^2-I_5h^2 +  2I_1I_3hg -  2I_1I_2I_3hg-2I_1I_3h^2-I_2^2h+I_3h^2\nonumber\\
&  \quad -  2I_3hg-I_4I_5hg-I_1I_2I_5h^2+I_2^2I_3g^2+I_1I_5h^2+I_4g,
\label{ch111}
\end{align}
\begin{align}
\varphi_2(r) & = 1 + I_5h-3I_2I_5h+I_1I_3h+I_3g-3I_2I_3g+I_4I_5g-I_3h-I_1I_4^2I_5g\nonumber\\
&  \quad+   3I_2^2I_3g- 2I_2I_4I_5g+2I_2I_3h-I_2^3I_5h+I_1I_2^2I_3h-I_2^3I_3g+3I_2^2I_5h\nonumber\\
&  \quad +  I_4I_5h-I_2I_4I_5h+2I_1I_3I_4h+I_3I_4g-I_2I_3I_4g+I_4^2+  2I_2^2I_4+I_2^4 \nonumber\\
&  \quad+  I_1I_2I_4I_5h-I_1^2I_3I_4h-I_1I_3I_4g + I_1I_2I_3I_4g- I_2^2I_3h -2I_1I_2^2I_4\nonumber\\
&  \quad-   2I_1I_4^2+I_1^2I_4^2+2I_4- 2I_1I_4-4I_2I_4+4I_1I_2I_4-4I_2 +  I_4^2I_5g \nonumber\\
&  \quad+   I_2^2I_4I_5g - 2I_1I_2I_3h-I_3I_4h+6I_2^2-  4I_2^3- I_1I_4I_5h
\label{ch112}
\end{align}
and
\begin{align}
\psi_2(r) & = h+I_3hg-2I_2I_4g+I_5h^2-3I_2h-2I_2I_3hg-I_1I_2I_3h^2+I_4g+  3I_2^2h \nonumber\\
&  \quad +  I_2^2I_4g+I_2^2I_5h^2-I_2^3h-I_3I_4hg+I_4^2g+I_4h -  I_2I_4h-2I_2I_5h^2\nonumber\\
&  \quad -  I_1I_4h+I_1I_2I_4h+I_1I_3h^2 + I_3I_4g^2+ I_2I_3h^2-I_3h^2+I_1I_3I_4hg\nonumber\\
&  \quad -  I_2I_3I_4g^2-2I_2I_4I_5hg +2I_4I_5hg  +I_2^2I_3hg-I_1I_4^2g+I_4^2I_5g^2.
\label{ch113}
\end{align}
The following Theorem is the main result of this work.
\begin{theorem}
\label{ch1T1}
Suppose that $h(q)$ and $g(q)$ are continuous functions such that
$\displaystyle\int_0^1 g(q)^2 \textrm{d}q = 1 $ and $\varphi_2(r)$ is a strictly
positive function on $[0,1)$. Then the equality
\begin{equation}
\label{ch114}
L[U](\nu) = U(\nu)+\displaystyle\int_0^{\nu} \int_0^r
\frac{\varphi_1(r)\;h(q)+
  \psi_2(r)\;g(q)}{\varphi_2(r)}\;\textrm{d}U(q)\;\textrm{d}r=w_{\nu}
\end{equation}
holds. Here $w_{\nu}, \;0\leq \nu \leq 1$ is a Wiener process.
\end{theorem}

\begin{pf} The proof will be done in
several steps.

\textbf{\textit{Step 1:}}  Introduce a Gaussian process
$$ M_t = \int_0^t q(t,s) \;\textrm{d}U(s),\qquad 0\leq t \leq 1,
$$
where the function $q(t,s)$ is chosen as solution of Fredholm equation
described in the next step. Observe that
$$
M_t = \int_0^t q(t,u) \;\textrm{d}W(u) - \int_0^1 g(u) \;\textrm{d}W(u)\int_0^t q(t,u) \; h(u)\;\textrm{d}u.
$$

\textbf{\textit{Step 2:}}  For the  correlation function of $M_t$
$$
R(t,s) = \textbf{\textrm{E}}\left[M_t M_s\right] ,\qquad t>s,
$$
we have
\begin{align*}
\textbf{\textrm{E}}\left[M_t M_s\right] & =  \textbf{\textrm{E}}\left[ \int_0^t q(t,u) \;\textrm{d}W(u) - \int_0^1  g(u) \;\textrm{d}W(u)\int_0^t q(t,u)\;h(u)\;\textrm{d}u\right]\\
&  \qquad \;\left[\int_0^s q(s,v) \;\textrm{d}W(v) - \int_0^1 g(v) \;\textrm{d}W(v)\int_0^s q(s,v) \;h(v)\;\textrm{d}v\right]\\
& =  \int_0^s q(t,u)\;q(s,u)\;\textrm{d}u - \int_0^t q(t,u)\; g(u)\;\textrm{d}u \int_0^s q(s,v)\; h(v)\;\textrm{d}v \\
&  \quad\quad- \; \int_0^s q(s,v) \;g(v)\;\textrm{d}v \int_0^t q(t,u) \;h(u)\;\textrm{d}u  \\
&  \quad\quad +  \underbrace{\int_0^1 g(v)^2 \textrm{d}v}_{=1} \int_0^t q(t,u)
\;h(u)\;\textrm{d}u \int_0^s q(s,v) \;h(v)\;\textrm{d}v\\
& =  \int_0^s q(s,u)\biggl[q(t,u) - \int_0^t q(t,v)\; g(v)\;\textrm{d}v\;h(u) \\
& \quad\quad- \int_0^t q(t,v)\; h(v)\;\textrm{d}v\; g(u)+ \int_0^t q(t,v)\; h(v)\;\textrm{d}v\; h(u)\biggl]\textrm{d}u\\
& =  \int_0^s q(s,u)\biggl[q(t,u) - \int_0^t q(t,v)\;\bigl[ g(v)\;h(u)\\
&  \qquad  \qquad \qquad\qquad + \;h(v)\;g(u)-h(v)\;h(u)\bigl]\;\textrm{d}v \biggl]\textrm{d}u.
\end{align*}
 Denote the  kernel
$$
K(u,v) = g(v)\;h(u)+h(v)\;g(u)-h(v)\;h(u).
$$
Then
\begin{equation}
\label{ch115}
\textbf{\textrm{E}}\left[M_t M_s\right] = \int_0^s q(s,u)\biggl[q(t,u) - \int_0^t q(t,v)\;K(u,v)\;\textrm{d}v \biggl]\textrm{d}u.
\end{equation}
Therefore if we take $q(t,s)$ such that it solves the Fredholm equation of the second kind
$(t$ is fixed$)$
\begin{equation}
\label{ch116} q(t,s) - \int_0^t q(t,v)\;K(s,v)\;\textrm{d}v = 1,
\qquad s \in[0,t],
\end{equation}
then $(\ref{ch115})$ becomes
\begin{equation}
\label{ch117} \textbf{\textrm{E}}\left[M_t M_s\right]  = \textbf{\textrm{E}}\left[M_s^2\right] = \int_0^s
q(s,u)\;\textrm{d}u.
\end{equation}

\textbf{\textit{Step 3:}}
The solution $q(t,s)$ can be found as follows. We have
$$
q(t,s) = 1+ \int_0^t q(t,v)\;K(s,v)\;\textrm{d}v.
$$
Denote
$$
A(t) = \int_0^t q(t,v)\; h(v)\;\textrm{d}v
$$
and
$$
B(t) = \int_0^t q(t,v)\; g(v)\;\textrm{d}v.
$$
Then $q(t,s)$ has the representation
\begin{equation}
\label{ch18}
q(t,s) = 1+ B(t) \; h(s) + A(t) \left( g(s)-
h(s)\right),
\end{equation}
where the function $A(t)$ itself is solution of the
following equation (after multiplying \eqref{ch18} by $h(s)$ and
integrating)
\begin{align}
\begin{split}
\label{ch118}
\int_0^t h(s) \; \textrm{d}s &= A(t) - B(t) \int_0^t h(s)^2\textrm{d}s \\
&  -  A(t) \int_0^t h(s)\; g(s)\;\textrm{d}s + A(t) \int_0^t h(s)^2\textrm{d}s.
\end{split}
\end{align}
 The function
$B(t)$ is solution of the following equation (after multiplying
\eqref{ch18} by $g(s)$ and integrating)
\begin{align}
\begin{split}
\label{ch119}
\int_0^t g(s) \; \textrm{d}s & = B(t) - B(t) \int_0^t h(s)\; g(s)\;\textrm{d}s \\
& - \;A(t) \int_0^t g(s)^2\textrm{d}s + A(t) \int_0^t h(s)\; g(s)\;\textrm{d}s.
\end{split}
\end{align}
Using the notation  \eqref{ch110},
we can write \eqref{ch118}$-$\eqref{ch119} as follows:
\begin{equation}
\label{ch120} A(t) - B(t)\; I_4 - A(t) \;I_2 + A(t)\; I_4 = I_3
\end{equation}
and
\begin{equation}
\label{ch121} B(t) - B(t)\; I_2 - A(t)\; I_1 + A(t) \; I_2 = I_5.
\end{equation}
Further, we have to find the expressions of  $A(t)$ and $B(t).$
Therefore we obtain from \eqref{ch121}
\begin{equation}
\label{ch122}
B(t) = \displaystyle \frac{I_5 + A(t) \left(I_1 -
I_2 \right)}{\left(1
-I_2\right)}.
\end{equation}
Then we insert \eqref{ch122} in \eqref{ch120} and obtain
$$ A(t) = \displaystyle \frac{I_3 \left(1- I_2\right)+I_4\;I_5}{\left(1-
I_2\right)^2+I_4\;I_6}
$$
and
$$ B(t) = \displaystyle \frac{ I_5\left(1+I_4-I_2\right)+ I_3\left(I_1 - I_2\right)}{\left(1-
I_2\right)^2+I_4\;I_6},
$$
where
$$
I_6 =\displaystyle\int_t^1 g(s)^2\;\textrm{d}s.
$$
Therefore the solution $q(t,s)$ of \eqref{ch116} is
\begin{align*}
q(t,s)&= 1+\displaystyle \frac{ I_5\left(1+I_4-I_2\right)+ I_3\left(I_1 - I_2\right)}{\left(1-
I_2\right)^2+I_4\;I_6}\; h(s)\\
 & \qquad + \displaystyle \frac{I_3 \left(1- I_2\right)+I_4\;I_5}{\left(1-
I_2\right)^2+I_4\;I_6}\;\left(g(s) - h(s)\right).
\end{align*}

 The final expression of $q(t,s)$ is
\begin{align*}
&q(t,s)= 1+\displaystyle \frac{
\displaystyle\int_0^tg(s)\textrm{d}s\left(1+\displaystyle\int_0^th(s)^2\textrm{d}s-\displaystyle
\int_0^th(s)g(s)\;\textrm{d}s\right)h(s)}{\left(1-\displaystyle\int_0^th(s)g(s)\;\textrm{d}s\right)^2
+\displaystyle\int_0^th(s)^2\textrm{d}s\displaystyle\int_t^1g(s)^2\textrm{d}s} \\
&  \qquad\qquad\;\;+ \displaystyle\frac{\displaystyle\int_0^th(s)\textrm{d}s\left(\displaystyle\int_0^tg(s)^2\textrm{d}s -\displaystyle\int_0^th(s)g(s)\;\textrm{d}s\right)h(s)}{\left(1-\displaystyle\int_0^th(s)g(s)\;\textrm{d}s\right)^2
+\displaystyle\int_0^th(s)^2\textrm{d}s\displaystyle\int_t^1g(s)^2\textrm{d}s} \\
& \quad+ \displaystyle \frac{\displaystyle\int_0^th(s)\textrm{d}s \left(1- \displaystyle\int_0^th(s)g(s)\textrm{d}s\right)+\displaystyle\int_0^th(s)^2\textrm{d}s\int_0^tg(s)\textrm{d}s}
{\left(1-\displaystyle\int_0^th(s)g(s)\textrm{d}s\right)^2
+\displaystyle\int_0^th(s)^2\textrm{d}s\displaystyle\int_t^1g(s)^2\textrm{d}s}(g(s) -h(s))
\end{align*}

\textbf{\textit{Step 4:}} To show that $M_t$ is martingale we need the
following Lemma.
\begin{Lem} \label{ch1L2} We have the following equality
\begin{equation}
\label{ch123}
\displaystyle\int_0^t q(t,s)\;\textrm{d}s = \displaystyle\int_0^t q(s,s)^2\textrm{d}s, \qquad 0\leq t \leq 1.
\end{equation}
\end{Lem}

\begin{pf}
Show that
$$
\frac{\textrm{d}}{\textrm{d}t}\int_0^t q(t,s)\;\textrm{d}s = \frac{\textrm{d}}{\textrm{d}t}\int_0^t q(s,s)^2\textrm{d}s = q(t,t)^2.
$$
Denote
\begin{align*}
C(t) & =  1 - \displaystyle\int_0^t h(s)\; g(s)\;\textrm{d}s,\\
D(t) & =
\displaystyle\int_0^tg(s)\;\textrm{d}s\left(1+\displaystyle\int_0^th(s)^2\textrm{d}s
-\displaystyle\int_0^th(s)\;g(s)\;\textrm{d}s\right)\\
&  \quad +  \displaystyle\int_0^th(s)\;\textrm{d}s\left(\displaystyle\int_0^tg(s)^2\textrm{d}s -
\displaystyle\int_0^th(s)\;g(s)\;\textrm{d}s\right),\\
K(t) & =  C(t)^2 + \displaystyle\int_0^t
h(s)^2\textrm{d}s\displaystyle\int_t^1 g(s)^2\textrm{d}s \\
\textrm{and}\qquad\qquad\qquad\;\\
N(t) & =  C(t)\; \int_0^t h(s)\; \textrm{d}s + \displaystyle\int_0^t
h(s)^2\textrm{d}s\displaystyle\int_0^t g(s)\;\textrm{d}s.
\end{align*}
Then $q(t,s)$ has the following expression
$$
q(t,s) = 1+\frac{D(t)}{K(t)}\;h(s) + \frac{N(t)}{K(t)}\;(g(s) - h(s)).
$$
Hence
\begin{align}
\begin{split}
\frac{\textrm{d}}{\textrm{d}t}\displaystyle\int_0^t q(t,s)\;\textrm{d}s & =  1+
\frac{D(t)\;h(t) +N(t)\;(g(t) - h(t))}{K(t)} \\
&  \quad+\frac{\Big(D'(t)\;K(t) -D(t)\;K'(t)\Big)}{K(t)^2}\int_0^th(s)\;\textrm{d}s\\
&  \quad +  \frac{\big(N'(t)K(t) -
N(t)K'(t)\big)}{K(t)^2}\int_0^t(g(s) - h(s))\;\textrm{d}s.\label{ch124}
\end{split}
\end{align}
Then we obtain the equalities
\begin{eqnarray*}
C(t) & =&  1 - I_2,\\
D(t) & =&  I_5\left(1+I_4-I_2\right)+ I_3\left(I_1 - I_2\right),\\
K(t) & =& \left(1 - I_2\right)^2 + I_4 \left(1-I_1\right),\\
N(t) & =&  I_3 \left(1 - I_2\right)+ I_4 \;I_5.
\end{eqnarray*}
The derivatives of these functions w.r.t. $t$ have such expressions
\begin{eqnarray*}
C'(t) & = &  - h(t)\; g(t),\\
D'(t) & = &g(t) + I_4 \; g(t) - I_2\; g(t) + I_5 \;h(t)^2 - I_5 \;h(t)\;g(t) \\
& & \qquad+ \;  I_1 \;h(t) - I_2 \;h(t)+I_3 \;g(t)^2  - I_3\;h(t)\;g(t),\\
K'(t) & = & -2 h(t) \;g(t) + 2I_2 \;h(t) g(t) + h(t)^2 - I_1 \;h(t)^2 - I_4 \;g(t)^2,\\
N'(t) & = & -I_3\;h(t)\; g(t) + I_4\;g(t)  + I_5 \;h(t)^2 + h(t) - I_2 \;h(t).
\end{eqnarray*}
Thus, \eqref{ch124} has the following representation
\begin{equation}
\label{ch125}
\frac{\textrm{d}}{\textrm{d}t}\displaystyle\int_0^t q(t,s)\;\textrm{d}s = 1+\frac{\Phi_1(t)+\Phi_2(t)+\Phi_3(t)}{K(t)^2},
\end{equation}
where
\begin{align*}
\Phi_1(t)  &=   \big(D(t)\;h(t) +N(t)\;(g(t) - h(t))\big)K(t),\\
\Phi_2(t) & = \big(D'(t)\;K(t) -D(t)\;K'(t)\big)I_3
\end{align*}
and
$$
\Phi_3(t) =  \big(N'(t)\;K(t) - N(t)\;K'(t)\big)\left(I_5- I_3\right).
$$
Then returning to initial notation  we obtain the
expressions
\begin{align*}
\Phi_1(t) & =   3I_2^2I_5h-3I_2I_5h+I_1I_3h - 3I_2I_3g+ I_1I_2^2I_3h-I_2^3I_5h + I_3g - I_2^3I_3g\\
&   \; + I_4I_5g-I_1I_3I_4g-I_3h-I_1^2I_3I_4h +2I_2I_3h+ I_5h- I_1I_4^2I_5g- I_2I_3I_4g  \\
&   \; -  2I_1I_2I_3h+3I_2^2I_3g- 2I_2I_4I_5g + I_4I_5h-I_2I_4I_5h+I_1I_2I_3I_4g+I_4^2I_5g\\
&   \; +  I_2^2I_4I_5g -  I_2^2I_3h+ 2I_1I_3I_4h-  I_3I_4h - I_1I_4I_5h + I_3I_4g+  I_1I_2I_4I_5h,  \\
\Phi_2(t) & =  I_2I_3^2h^2+ 2I_3I_4g-3I_2I_3g + I_3I_5hg+I_3I_4^2I_5g^2+I_3^2I_4g^2+ I_1I_3h\\
&  \; -  3I_2I_3I_4g+3I_2^2I_3g-I_2I_3I_5h^2-2I_2I_3I_5hg+I_3^2g^2-I_1I_3^2h^2\\
& \; +  I_2^2I_3I_4g-I_2^3I_3g+I_2^2I_3I_5h^2+I_2^2I_3I_5hg+I_1I_2^2I_3h-2I_1I_2I_3h\\
&  \; -  I_3^2hg-2I_2I_3^2g^2+I_2^2I_3^2g^2+2I_2^2I_3h-2I_2I_3I_4I_5hg-I_1^2I_3I_4h\\
&  \;+  I_2^2I_3^2hg+I_3I_4^2g+I_3I_4I_5hg-2I_1I_2I_3^2hg+I_1^2I_3^2h^2- I_1I_2I_3^2h^2\\
&  \; -  I_1I_3I_4g-I_1I_3I_4^2g+I_1I_2I_3I_4g+I_1I_3I_4I_5hg-I_1I_2I_3I_5h^2+I_3g\\
&  \; + I_1I_3^2I_4hg + I_1I_3I_5h^2+I_3I_4I_5g^2+I_1I_3I_4h-I_2I_3h-I_2I_3^2I_4g^2 \\
& \; -  I_2I_3I_4I_5g^2 + 2I_1I_3^2hg-I_3^2I_4hg+I_1I_2I_3I_4h-I_2I_3I_4h-I_2^3I_3h
\end{align*}
and
\begin{align*}
\Phi_3(t) & =  I_3I_5hg+  I_2I_3I_4h-I_3I_5h^2 -3I_2I_5h-I_3I_4^2g-I_2^2I_3^2hg\\
&  \;- 2 I_2I_5^2h^2 + 3I_2^2I_5h+I_2^2I_3I_5hg- I_2^2I_3I_5h^2+I_2^2I_5^2h^2-I_3I_4^2I_5g^2\\
&  \; +  I_4^2I_5g-3I_2^2I_3h -  I_2I_4I_5h+I_1I_3I_4I_5hg-  I_3^2hg -I_1I_4^2I_5g\\
&  \; +I_2I_3^2I_4g^2  +I_1I_3I_5h^2+I_3I_4I_5g^2 + I_2I_3I_5h^2+ I_5h-2I_2I_3I_5hg  \\
&  \; +  2I_4I_5^2hg -  2I_2I_4I_5^2hg+I_4^2I_5^2g^2+I_2^2I_4I_5g-I_3h+  I_1I_2I_4I_5h  \\
&  \; + 3I_2I_3h+2I_2I_3^2hg+2I_2I_3I_4g +  2I_2I_3I_5h^2+I_5^2h^2-  I_1I_4I_5h \\
&  \; -  I_2^2I_3I_4g +I_2^3I_3h +  I_3^2I_4hg+I_4I_5h-I_2^3I_5h+I_3^2h^2+I_4I_5g\\
&  \; -  I_1I_3^2I_4hg+  I_1I_2I_3^2h^2 +  I_1I_3I_4^2g+I_1I_3I_4h-I_3I_4h-I_3I_4g \\
&  \; -  I_3^2I_4g^2 -  I_2I_3^2h^2 -I_3I_5h^2-I_1I_2I_3I_4h -I_2I_3I_4I_5g^2-  I_1I_3^2h^2\\
&  \; -  2I_3I_4I_5hg +  2I_2I_3I_4I_5hg-2I_2I_4I_5g-  I_1I_2I_3I_5h^2-  I_3I_4I_5hg .
\end{align*}
 Denote $$\Phi(t) = \Phi_1(t)+\Phi_2(t)+\Phi_3(t),$$ then we can write
\begin{align*}
\Phi(t) & = h(t)\;\Big(2I_5-6I_2I_5+2I_1I_3-2I_3+6I_2^2I_5-4I_1I_2I_3+4I_2I_3+  4I_1I_3I_4\\
&  \;-2I_2^3I_5+2I_1I_2^2I_3+  2I_4I_5 -  2I_2I_4I_5  -  2I_3I_4+2I_1I_2I_4I_5-2I_1^2I_3I_4\\
&  \;-2I_1I_4I_5-2I_2^2I_3\Big)+\;g(t)\;\Big(2I_3- 2I_1I_3I_4+2I_1I_2I_3I_4-  6I_2I_3+2I_4I_5\\
&  \;-2I_1I_4^2I_5 +  6I_2^2I_3 - 4I_2I_4I_5-  2I_2^3I_3+2I_2^2I_4I_5+2I_3I_4+2I_4^2I_5 \\
&  \;-  2I_2I_3I_4\Big)+  h(t)g(t)
\Big(2I_3I_5 - 2I_3^2-4I_2I_3I_5-  2I_3I_4I_5 + 2I_1I_3I_4I_5\\
& \;+2I_2I_3^2+2I_4I_5^2+2I_1I_3^2-2I_1I_2I_3^2 - 2I_2I_4I_5^2+2I_2^2I_3I_5\Big)+ h(t)^2\\&\;\Big(2I_1I_3I_5- 2I_1I_2I_3I_5 -  2I_2I_5^2+I_2^2I_5^2-2I_3I_5 +  2I_2I_3I_5-2I_1I_3^2+I_3^2\\
&  \;+I_5^2+ I_1^2I_3^2\Big)+  g(t)^2
\Big(I_3^2-2I_2I_3^2+2I_3I_4I_5-2I_2I_3I_4I_5+I_4^2I_5^2+I_2^2I_3^2\Big).
\end{align*}
Finally, \eqref{ch125} can be written as follows:
\begin{equation}
\label{ch126}
\frac{\textrm{d}}{\textrm{d}t}\displaystyle\int_0^t q(t,s)\;\textrm{d}s = 1+\frac{\Phi(t)}{K(t)^2}.
\end{equation}
Now the expression of $q(t,t)^2$ is
\begin{eqnarray*}
q(t,t)^2& = &\left(1+\frac{D(t)}{K(t)}\;h(t) +
\frac{N(t)}{K(t)}\;(g(t) - h(t))\right)^2\\
& = & 1+ \frac{2\;h(t)\;K(t)\big(D(t) - N(t)\big)\; + 2 \;g(t)\; K(t)\;N(t)}{K(t)^2}\\
& & \quad+\; \frac{h(t)^2\big(D(t) -N(t)\big)^2+ g(t)^2 N(t)^2}{K(t)^2}\\
& & \quad+\;\frac{2\;h(t)\;g(t) N(t) \big(D(t) - N(t)\big)}{K(t)^2}.
\end{eqnarray*}
 Denote
\begin{eqnarray*}
M(t) & = & 2\;K(t)\left(D(t) - N(t)\right), \\
Q(t) & = & 2\; K(t)\;N(t),\\
L(t) & = & \big(D(t) -N(t)\big)^2, \\
E(t) & = & N(t)^2,\\
H(t) & = & 2\;N(t) \big(D(t)-N(t)\big).
\end{eqnarray*}
Therefore the final expression for $q(t,t)^2$ is
\begin{align}
\begin{split}
\label{ch127}
q(t,t)^2 & =  1+ \frac{h(t)\;M(t) +  g(t)\; Q(t)+h(t)^2L(t) }{K(t)^2}\\
&  \quad +\;\frac{g(t)^2 E(t)+ h(t)\;g(t) \;H(t)}{K(t)^2},
\end{split}
\end{align}
where
\begin{align*}
M(t) & = 2I_1I_3-6I_2I_5-2I_3+6I_2^2I_5-4I_1I_2I_3+4I_2I_3-2I_2^3I_5+  2I_4I_5\\
&  \; +2I_1I_2^2I_3-2I_1^2I_3I_4  -   2I_2^2I_3 -2I_2I_4I_5+2I_5-2I_3I_4-2I_1I_4I_5\\
&  \; + 2I_1I_2I_4I_5+4I_1I_3I_4,\\
Q(t) & =  2I_3I_4 -6I_2I_3 + 6I_2^2I_3 -2I_2^3I_3 + 2I_4I_5 - 4I_2I_4I_5 + 2I_2^2I_4I_5\\
&  \;  + 2I_3-2I_2I_3I_4- 2I_1I_3I_4+2I_1I_2I_3I_4+2I_4^2I_5-2I_1I_4^2I_5,\\
L(t) & =  I_5^2+I_2^2I_5^2-2I_2I_5^2+I_1^2I_3^2+I_3^2-2I_1I_3^2+2I_1I_3I_5-2I_3I_5\\
&  \; +  2I_2I_3I_5-2I_1I_2I_3I_5, \\
E(t) & =  I_3^2-2I_2I_3^2+I_2^2I_3^2+I_4^2I_5^2+2I_3I_4I_5-2I_2I_3I_4I_5,\\
H(t) & =  2I_3I_5-4I_2I_3I_5-2I_3I_4I_5+2I_1I_3^2+2I_2I_3^2+2I_2^2I_3I_5\\
& \; -  2I_2I_4I_5^2 +  2I_1I_3I_4I_5-2I_3^2-2I_1I_2I_3^2+2I_4I_5^2.
\end{align*}
The comparison of these expressions with \eqref{ch126}$-$\eqref{ch127} shows that
the Lemma is proved.
\end{pf}

\textbf{\textit{Step 5:}} In the next step, we need the following Lemma to show that the linear transformation is a Wiener process.
\begin{Lem}
\label{ch1L3}
If the Gaussian process $M_s$ satisfies \eqref{ch117} and we have relation \eqref{ch123}, then
$$z(t) = \displaystyle\int_0^t q(s,s)^{-1}\;dM_s$$
is a Wiener process.
\end{Lem}
\textbf{Proof.} The proof can be found, e.g., in \cite{KK14}, Lemma 2.

Hence
$$
M_t = \int_0^t q(s,s) \;\textrm{d}w_s,\qquad0\leq t\leq 1
$$
is a Gaussian martingale,
where $w_s,0\leq s \leq 1$ is a Wiener process.

Therefore we have the  equality
$$
\textit{w}_t = \displaystyle\int_0^t q(s,s)^{-1}\;\textrm{d}M_s =
U(t)+\displaystyle\int_0^t q(s,s)^{-1}\int_0^s q_s'(s,v)\;\textrm{d}U(v)\;\textrm{d}s,
$$
where $w_t, 0 \leq t \leq 1$ is a Wiener process (by Lemma 3).

Now we have to calculate the right side of the above expression. The derivative $q_t'(t,s)$ w.r.t. $t$ can be written as follows:
$$
q_t'(t,s) = \frac{\left(\psi_1(t) - \psi_2(t)\right)h(s)+ \psi_2(t)\; g(s)}{K(t)^2},
$$
where
$$
\psi_1(t) = D'(t)\;K(t) -D(t)\;K'(t)$$and
$$\psi_2(t)  =  N'(t)\;K(t) - N(t)\;K'(t).$$
Returning to initial notation we obtain the following expression:
\begin{align*}
\psi_1(t) & = g+2I_4g-3I_2g+I_5hg+I_1h-I_2h+I_1I_4I_5hg+I_4^2g+3I_2^2g\\
&  \quad -\;  I_2I_5h^2 -2I_2I_5hg-2I_1I_2h+2I_2^2h+I_2^2I_3g^2-I_1I_2I_3h^2-I_2^3g\\
&  \quad +\;  I_2^2I_5hg+I_1I_2^2h-I_2^3h-2I_2I_4I_5hg+I_1^2I_3h^2-I_1I_2I_5h^2-I_3hg\\
&  \quad -\;  I_2I_4h+I_3I_4g^2-I_3I_4hg-I_1I_4g-I_1I_4^2g+I_1I_2I_4g-2I_2I_3g^2\\
&  \quad +\;  I_1I_2I_4h+I_1I_5h^2+I_4I_5g^2+I_4^2g+I_2^2I_3hg-I_1I_3h^2+I_4^2I_5g^2\\
&  \quad -\;  I_2I_4I_5g^2+2I_1I_3hg-2I_1I_2I_3hg+I_2^2I_4g-I_1^2I_4h+I_2I_3h^2\\
&  \quad -\;  I_2I_3I_4g^2-3I_2I_4g+I_2^2I_5h^2+I_1I_4h+I_1I_3I_4hg+I_4I_5hg+I_3g^2.
\end{align*}
The function  $\psi_2(t)$ is defined by \eqref{ch113}.
Hence we obtain
$$
\frac{q_t'(t,s)}{q(t,t)} = \frac{\left(\psi_1(t) - \psi_2(t)\right)\;h(s)+ \psi_2(t)\; g(s)}{K(t)^2+ \Phi_1(t)}.
$$
Then if we put
$$
\varphi_1(t)  = \psi_1(t) - \psi_2(t)$$ and
$$\varphi_2(t)  =  K(t)^2+ \Phi_1(t),$$
then this implies that $$\frac{q_t'(t,s)}{q(t,t)} = \frac{\varphi_1(t)\;h(s)+ \psi_2(t)\; g(s)}{\varphi_2(t)},$$
with $\varphi_1(t)$ and $\varphi_2(t)$ were defined by \eqref{ch111}$-$\eqref{ch112}.
Finally, we obtain the  expression
$$
w_t = \displaystyle\int_0^t q(s,s)^{-1}\textrm{d}M_s =
U(t)+\displaystyle\int_0^t \int_0^s \frac{\varphi_1(s)\;h(v)+ \psi_2(s)\;
 g(v)}{\varphi_2(s)}\;\textrm{d}U(v)\;\textrm{d}s.
$$
 This is the linear transformation
$w_t = L[U](t)$ of the process $U(\cdot)$ into the Wiener process
$w_t$ and this proves the Theorem \ref{ch1T1}.
\end{pf}

\textbf{Remark.} Let us present a sufficient condition for $\varphi_2(t)>0$.

$\mathcal{R}_0.$  \textit{Suppose that $h(t)$ and $g(t)$ are continuous strictly positive functions such that
$
g(t)>h(t)
$
 and
$$
\int_0^th(s)\;g(s)\;\textrm{d}s <1,\qquad \int_0^1 g(s)^2 \textrm{d}s = 1, \qquad  0\leq t <1,
$$
then
$\varphi_2(t)$ defined by \eqref{ch112} is strictly positive function on $[0,1)$.}

Now we will verify that if the condition $\mathcal{R}_0$ is satisfied, then $\varphi_2(t)$ is strictly positive function. Remind that $\varphi_2(t)$ has the following expression:
$$
\varphi_2(t)  =  K(t)^2+ \Phi_1(t),\qquad \qquad  0\leq t <1,
$$
where
$$
\Phi_1(t)  =   \big(D(t)\;h(t) +N(t)\;(g(t) - h(t))\big)K(t),
$$
with
\begin{align*}
D(t) & =
\displaystyle\int_0^tg(s)\;\textrm{d}s\left(1+\displaystyle\int_0^th(s)^2\textrm{d}s
-\displaystyle\int_0^th(s)\;g(s)\;\textrm{d}s\right)\\
&  \quad +  \displaystyle\int_0^th(s)\;\textrm{d}s\left(\displaystyle\int_0^tg(s)^2\textrm{d}s -
\displaystyle\int_0^th(s)\;g(s)\;\textrm{d}s\right),
\end{align*}
$$
N(t)  =  C(t)\; \int_0^t h(s)\; \textrm{d}s + \displaystyle\int_0^t
h(s)^2\textrm{d}s\displaystyle\int_0^t g(s)\;\textrm{d}s
$$
and
$$
K(t)  =  C(t)^2 + \displaystyle\int_0^t
h(s)^2\textrm{d}s\displaystyle\int_t^1 g(s)^2\textrm{d}s.
$$
Here
$$
C(t)=  1 - \displaystyle\int_0^t h(s)\; g(s)\;\textrm{d}s.
$$
Note that it is sufficient to check that $\Phi_1(t) >0$ to obtain $\varphi_2(t)>0$. Recall that $K(t)$ is strictly positive function. Then due to the following condition
\begin{equation}
\label{ch16}
\displaystyle\int_0^th(s)\;g(s)\;\textrm{d}s <1,
\end{equation}
we have $C(t)>0$. Consequently, we obtain $N(t)>0$ by the conditions
\begin{equation}
\label{ch17}
 h(t)>0 \qquad \textrm{and} \qquad g(t)>0.
\end{equation}
Finally, we see that $D(t)>0$ by the conditions $\displaystyle\int_0^1 g(s)^2 \textrm{d}s = 1$, \eqref{ch16} and \eqref{ch17}.
We conclude that we have $\varphi_2(t)>0$ if we suppose that $g(t)>h(t)$.

\section{Test}
\label{ch1s04}

Our objective is to test the composite parametric hypothesis ${\cal H}_0$ and
to do this we will propose a statistic based on the MDE
$\vartheta _\varepsilon ^*$. Recall that the starting statistic
\begin{equation}
\label{ch1010}
u_\varepsilon \left(t\right)=\frac{X_t-x_t\left(\vartheta _\varepsilon
  ^*\right)}{\varepsilon \;S\left(\vartheta _\varepsilon
  ^*,X_t \right)}
 \end{equation}
  converges to the random function
 \begin{align*}
u\left(t\right)=\int_{0}^{t}\frac{{\rm d}W_s}{S(\vartheta,x_s)}-\int_0^T
\int_v^T\frac{S(\vartheta,x_s)\;\dot{x}_s(\vartheta)}{ J(\vartheta)  S\left(\vartheta
  ,x_v\right)   }\;\textrm{d}s\;\textrm{d}W_v\int_0^t
\frac{\dot{S}(\vartheta,x_s)}{S(\vartheta,x_s)}\;\textrm{d}s .
 \end{align*}
Then the linear transformation
\begin{align*}
U\left(\frac{t}{T}\right)=\frac{1}{\sqrt{T}}\int_{0}^{t}S\left(\vartheta
,x_s\right){\rm d}u\left(s\right)
\end{align*}
has the following representation, by It\^{o} formula,
\begin{align*}
&U\left(\frac{t}{T}\right)  =  \frac{1}{\sqrt{T}}\;S(\vartheta,x_t)\;u(t)
- \frac{1}{\sqrt{T}}
\int_0^tS'(\vartheta,x_s)\;S(\vartheta,x_s)\;u(s)\;\textrm{d}s
\end{align*}
which leads to the random function
$$
U\left(\nu \right)
=W\left(\nu \right)-\int_{0}^{1}g\left(\vartheta ,r\right){\rm d}W\left(r\right)\int_{0}^{\nu }h\left(\vartheta ,r\right){\rm d}r.
$$
The last step is to apply the transformation $L\left[U\right](\cdot)$ from Theorem
\ref{ch1T1} and to obtain the Wiener process
$$
L\left[U\right]\left(\nu \right)=w_\nu ,\qquad 0\leq \nu \leq 1.
$$
Now we have to realize the similar transformations with the ``empirical''
process $u_\varepsilon \left(\cdot\right)$ defined by \eqref{ch1010}, i.e., we (formally) calculate
\begin{align}
\label{ch1000}
U_\varepsilon
\left(\frac{t}{T}\right)=\frac{1}{\sqrt{T}}\int_{0}^{t}S\left(\vartheta_\varepsilon
^*,X_s \right)\;{\rm d}u_\varepsilon \left(s\right).
\end{align}
Then we apply the transformation $L\left[\cdot \right]$ to the process
$U_\varepsilon\left(\cdot\right) $ and  we show that this statistic converges in distribution to the Wiener process $w_\nu , 0\leq \nu \leq 1$.
Therefore the test $\psi _\varepsilon =\1_{\left\{\Delta _\varepsilon>c_\alpha
  \right\}}$ with
$$
\Delta _\varepsilon =\frac{1}{T}\int_{0}^{T}L\left[U_\varepsilon \right]\left(t \right)^2{\rm d}t \Longrightarrow \int_{0}^{1}w_\nu ^2\;{\rm d}\nu
$$
will be ADF because the limit distribution of $\Delta _\varepsilon $ does not
depend on $S\left(\cdot,\cdot \right)$ and $\vartheta $.

Let us realize this program. We have the following representation for the process $U_{\varepsilon}(\cdot)$
\begin{align}
\begin{split}
\label{ch128}
&U_{\varepsilon}\left(\frac{t}{T}\right) =  \frac{1}{\sqrt{T}}\;S\left(\vartheta^{*}_{\varepsilon},X_t\right)\;u_\varepsilon\left(t\right)\\
&\qquad \qquad\qquad- \frac{1}{\sqrt{T}}
\int_0^tS'\left(\vartheta^{*}_{\varepsilon},X_s\right)\;S\left(\vartheta^{*}_{\varepsilon},X_s\right)
\;u_\varepsilon\left(s\right)\;\textrm{d}s.
\end{split}
\end{align}

Introduce the functions
\begin{equation}
\label{cch05}
\hat{h}\left(\vartheta,v\right)=\frac{1}{\sqrt{T}}\;J(\vartheta)^{-1}\dot{S}(\vartheta,x_v), \qquad J(\vartheta)=\int_0^T\dot{x}_s(\vartheta)^2\textrm{d}s
\end{equation}
 and
\begin{equation}
\label{cch06}
\hat{g}\left(\vartheta,v\right)=S(\vartheta,x_v)^{-1} I\left(\vartheta,v\right)
, \qquad I\left(\vartheta,v\right)=\int_v^TS(\vartheta,x_s)\;\dot{x}_{s}(\vartheta)\;\textrm{d}s
\end{equation}
and their ``empirical versions'', respectively
\begin{align*}
& h_{\varepsilon}\left(\vartheta^{*}_{\varepsilon},v\right) =
\frac{1}{\sqrt{T}}\;J_\varepsilon(\vartheta^{*}_{\varepsilon})^{-1}
\;\dot{S}(\vartheta^*_{\varepsilon},X_v),\qquad
g_{\varepsilon}\left(\vartheta^{*}_{\varepsilon},v\right)
=S(\vartheta^{*}_{\varepsilon},X_v)^{-1}I_{\varepsilon}(\vartheta^{*}_{\varepsilon},v).
\end{align*}
Here
\begin{align*}
&J_\varepsilon(\vartheta^{*}_{\varepsilon})=
\int_0^T\dot{x}_s(\vartheta^{*}_{\varepsilon})^2\textrm{d}s,\qquad
  I_{\varepsilon}(\vartheta^{*}_{\varepsilon},v)=\int_v^TS
(\vartheta^{*}_{\varepsilon},X_s)\;\dot{x}_s(\vartheta^{*}_{\varepsilon})\;\textrm{d}s.
\end{align*}
Note that in the functions $\hat{h}\left(\cdot,\cdot\right)$ and $\hat{g}\left(\cdot,\cdot\right)$ we omit the normalizing constants in the expressions of the functions $h\left(\cdot,\cdot\right)$ and $g\left(\cdot,\cdot\right)$ defined by \eqref{ch107}$-$\eqref{ch108}, for simplicity of exposition, because the structure of the used statistic is such that we can do this without changing the limit distribution of the statistic.

Then denote the ``empirical versions''
$$
I_{1,\varepsilon} =C_1\left(T\right) \int_0^s
g_{\varepsilon}\left(\vartheta^*_{\varepsilon},v\right)^2\textrm{d}v,
\quad I_{2,\varepsilon} =C_2\left(T\right) \int_0^s
h_{\varepsilon}\left(\vartheta^*_{\varepsilon},v\right)
g_{\varepsilon}\left(\vartheta^*_{\varepsilon},v\right)\;\textrm{d}v,
$$$$
I_{3,\varepsilon} =C_3\left(T\right)\int_0^s
h_{\varepsilon}\left(\vartheta^*_{\varepsilon},v\right)\;\textrm{d}v,
\quad I_{4,\varepsilon} =C_4\left(T\right)\int_0^s
h_{\varepsilon}\left(\vartheta^*_{\varepsilon},v\right)^2\textrm{d}v$$
and
$$I_{5,\varepsilon} =C_5\left(T\right)\int_0^s
g_{\varepsilon}\left(\vartheta^*_{\varepsilon},v\right)\;\textrm{d}v
$$
of the integrals
$$
\hat{I}_1 = C_1\left(T\right)\int_0^s
\hat{g}\left(\vartheta,v\right)^2\textrm{d}v,
\quad \hat{I}_2 = C_2\left(T\right)\int_0^s
\hat{h}\left(\vartheta,v\right)
\hat{g}\left(\vartheta,v\right)\;\textrm{d}v,
$$$$
\hat{I}_3 =C_3\left(T\right)\int_0^s
\hat{h}\left(\vartheta,v\right)\;\textrm{d}v,
\quad \hat{I}_4 =C_4\left(T\right)\int_0^s
\hat{h}\left(\vartheta,v\right)^2\textrm{d}v$$
and
$$\hat{I}_5 =C_5\left(T\right)\int_0^s
\hat{g}\left(\vartheta,v\right)\;\textrm{d}v,
$$
where $C_1\left(T\right)=\displaystyle\frac{1}{T^3}$, $C_2\left(T\right)=\sqrt{T}$, $C_3\left(T\right)=T\sqrt{T}$, $C_4\left(T\right)=T^4$ and $C_5\left(T\right)=\displaystyle\frac{1}{T^2}$.
This allows us to introduce the ``empirical versions''
$\varphi_{1,\varepsilon}(\cdot), \varphi_{2,\varepsilon}(\cdot)$ and
$\psi_{2,\varepsilon}(\cdot)$ of $\hat{\varphi}_{1}(\cdot), \hat{\varphi}_{2}(\cdot)$ and
$\hat{\psi}_2(\cdot)$ defined respectively by \eqref{ch111}$-$\eqref{ch113}, where we replace the functions $g$ by $\hat{g}=\displaystyle\;\frac{1}{T}\;\hat{g}\left(\vartheta,s\right)$ and $h$ by $\hat{h}=T^2\sqrt{T}\;\hat{h}\left(\vartheta,s\right)$,
\begin{align}
\label{ch10}
\varphi_{1,\varepsilon}(s) & =
g_{\varepsilon}-h_{\varepsilon}+I_{4,\varepsilon}g_{\varepsilon}-3I_{2,\varepsilon}g_{\varepsilon}
+I_{5,\varepsilon}h_{\varepsilon}g_{\varepsilon}+I_{1,\varepsilon}h_{\varepsilon}+2I_{2,\varepsilon}
h_{\varepsilon}
+I_{2,\varepsilon}^2I_{3,\varepsilon}g_{\varepsilon}^2 \nonumber\\ &\quad
-2I_{2,\varepsilon}I_{3,\varepsilon}g_{\varepsilon}^2-
I_{5,\varepsilon}h_{\varepsilon}^2 + I_{3,\varepsilon}g_{\varepsilon}^2-
I_{2,\varepsilon}I_{4,\varepsilon}g_{\varepsilon}+3I_{2,\varepsilon}^2g_{\varepsilon}
-2I_{3,\varepsilon}h_{\varepsilon}g_{\varepsilon}\nonumber\\ & \quad
-2I_{2,\varepsilon}I_{5,\varepsilon}h_{\varepsilon}g_{\varepsilon}
-I_{2,\varepsilon}^3g_{\varepsilon}+I_{2,\varepsilon}^2I_{5,\varepsilon}h_{\varepsilon}g_{\varepsilon}
+I_{1,\varepsilon}I_{2,\varepsilon}^2h_{\varepsilon}-I_{2,\varepsilon}^2h_{\varepsilon}
\nonumber\\ &\quad+2I_{1,\varepsilon}I_{4,\varepsilon}h_{\varepsilon}-I_{1,\varepsilon}I_{4,\varepsilon}g_{\varepsilon}
+ I_{1,\varepsilon}^2I_{3,\varepsilon}h_{\varepsilon}^2+
I_{1,\varepsilon}I_{4,\varepsilon}I_{5,\varepsilon}h_{\varepsilon}g_{\varepsilon}+
I_{2,\varepsilon}I_{5,\varepsilon}h_{\varepsilon}^2 \nonumber\\ & \quad
+I_{1,\varepsilon}I_{5,\varepsilon}h_{\varepsilon}^2
+I_{4,\varepsilon}I_{5,\varepsilon}g_{\varepsilon}^2-I_{1,\varepsilon}I_{2,\varepsilon}I_{5,\varepsilon}
h_{\varepsilon}^2-I_{4,\varepsilon}h_{\varepsilon}-
I_{2,\varepsilon}I_{4,\varepsilon}I_{5,\varepsilon}g_{\varepsilon}^2
\nonumber\\ & \quad
+2I_{1,\varepsilon}I_{3,\varepsilon}h_{\varepsilon}g_{\varepsilon} -
2I_{1,\varepsilon}I_{2,\varepsilon}I_{3,\varepsilon}h_{\varepsilon}g_{\varepsilon}-2I_{1,\varepsilon}
I_{3,\varepsilon}h_{\varepsilon}^2-2I_{1,\varepsilon}I_{2,\varepsilon}h_{\varepsilon}\nonumber\\ &
\quad +
I_{1,\varepsilon}I_{2,\varepsilon}I_{4,\varepsilon}g_{\varepsilon}+2I_{2,\varepsilon}I_{3,\varepsilon}
h_{\varepsilon}g_{\varepsilon}
-I_{4,\varepsilon}I_{5,\varepsilon}h_{\varepsilon} g_{\varepsilon}\nonumber
+I_{3,\varepsilon}h_{\varepsilon}^2-I_{1,\varepsilon}^2
I_{4,\varepsilon}h_{\varepsilon},\\
\varphi_{2,\varepsilon}(s) & = 1 -2I_{1,\varepsilon}I_{4,\varepsilon}^2
-3I_{2,\varepsilon}I_{5,\varepsilon}h_{\varepsilon}+I_{1,\varepsilon}I_{3,\varepsilon}h_{\varepsilon}
+I_{3,\varepsilon}g_{\varepsilon}-3I_{2,\varepsilon}I_{3,\varepsilon}g_{\varepsilon}
\nonumber\\ &
\quad+I_{4,\varepsilon}I_{5,\varepsilon}g_{\varepsilon} -
I_{3,\varepsilon}h_{\varepsilon}+3I_{2,\varepsilon}^2I_{5,\varepsilon}h_{\varepsilon}
-
2I_{1,\varepsilon}I_{2,\varepsilon}I_{3,\varepsilon}h_{\varepsilon}-2I_{1,\varepsilon}I_{2,\varepsilon}^2
I_{4,\varepsilon}
\nonumber\\ & \quad-
2I_{2,\varepsilon}I_{4,\varepsilon}I_{5,\varepsilon}g_{\varepsilon}+2I_{2,\varepsilon}I_{3,\varepsilon}
h_{\varepsilon}-
I_{2,\varepsilon}^3I_{5,\varepsilon}h_{\varepsilon}+I_{1,\varepsilon}I_{2,\varepsilon}^2I_{3,\varepsilon}
h_{\varepsilon}-
2I_{1,\varepsilon}I_{4,\varepsilon}\nonumber\\ &
\quad-I_{2,\varepsilon}^3I_{3,\varepsilon}g_{\varepsilon} +
I_{4,\varepsilon}^2 +
I_{2,\varepsilon}^2I_{4,\varepsilon}I_{5,\varepsilon}g_{\varepsilon}+I_{4,\varepsilon}I_{5,\varepsilon}
h_{\varepsilon}+
I_{4,\varepsilon}^2I_{5,\varepsilon}g_{\varepsilon}\nonumber\\ &\quad -
I_{2,\varepsilon}I_{4,\varepsilon}I_{5,\varepsilon}h_{\varepsilon}+2I_{1,\varepsilon}
I_{3,\varepsilon}I_{4,\varepsilon}h_{\varepsilon}+I_{3,\varepsilon}I_{4,\varepsilon}g_{\varepsilon}
+3I_{2,\varepsilon}^2I_{3,\varepsilon}g_{\varepsilon}
\nonumber\\ & \quad+
2I_{2,\varepsilon}^2I_{4,\varepsilon} + 2I_{4,\varepsilon}-
I_{3,\varepsilon}I_{4,\varepsilon}h_{\varepsilon}+6I_{2,\varepsilon}^2+
I_{1,\varepsilon}I_{2,\varepsilon}I_{4,\varepsilon}I_{5,\varepsilon}h_{\varepsilon}\nonumber\\ & \quad
-I_{1,\varepsilon}^2I_{3,\varepsilon}I_{4,\varepsilon}h_{\varepsilon}-I_{1,\varepsilon}I_{3,\varepsilon}
I_{4,\varepsilon}g_{\varepsilon}
+I_{5,\varepsilon}h_{\varepsilon}+
I_{1,\varepsilon}I_{2,\varepsilon}I_{3,\varepsilon}I_{4,\varepsilon}g_{\varepsilon}\nonumber\\ & \quad -
I_{1,\varepsilon}I_{4,\varepsilon}^2I_{5,\varepsilon}g_{\varepsilon}-
I_{1,\varepsilon}I_{4,\varepsilon}I_{5,\varepsilon}h_{\varepsilon}+4I_{1,\varepsilon}I_{2,\varepsilon}
I_{4,\varepsilon}
-I_{2,\varepsilon}I_{3,\varepsilon}I_{4,\varepsilon}g_{\varepsilon}\nonumber\\
&\quad -
I_{2,\varepsilon}^2I_{3,\varepsilon}h_{\varepsilon}+I_{2,\varepsilon}^4-
4I_{2,\varepsilon}^3-4I_{2,\varepsilon}-4I_{2,\varepsilon}I_{4,\varepsilon}+ I_{1,\varepsilon}^2I_{4,\varepsilon}^2
\end{align}
and
\begin{align*}
\psi_{2,\varepsilon}(s) & = h_{\varepsilon}+I_{3,\varepsilon}h_{\varepsilon}g_{\varepsilon}+I_{4,\varepsilon}g_{\varepsilon}
+I_{5,\varepsilon}h_{\varepsilon}^2-3I_{2,\varepsilon}h_{\varepsilon}-I_{1,\varepsilon}I_{2,\varepsilon}
I_{3,\varepsilon}h_{\varepsilon}^2
\nonumber\\
&  \quad-2I_{2,\varepsilon}I_{4,\varepsilon}g_{\varepsilon}+I_{4,\varepsilon}h_{\varepsilon} -2I_{2,\varepsilon}I_{5,\varepsilon}h_{\varepsilon}^2 +  3I_{2,\varepsilon}^2h_{\varepsilon}
+  I_{2,\varepsilon}^2I_{4,\varepsilon}g_{\varepsilon}\nonumber\\
&  \quad +I_{2,\varepsilon}^2I_{5,\varepsilon}h_{\varepsilon}^2-I_{2,\varepsilon}^3h_{\varepsilon}-I_{3,\varepsilon}
I_{4,\varepsilon}h_{\varepsilon}g_{\varepsilon}
+I_{4,\varepsilon}^2g_{\varepsilon} - I_{2,\varepsilon}I_{4,\varepsilon}h_{\varepsilon}\nonumber\\
&  \quad +I_{1,\varepsilon}I_{3,\varepsilon}I_{4,\varepsilon}h_{\varepsilon}g_{\varepsilon}-I_{1,\varepsilon}
I_{4,\varepsilon}^2g_{\varepsilon} -  I_{1,\varepsilon}I_{4,\varepsilon}h_{\varepsilon}+I_{1,\varepsilon}I_{2,\varepsilon}I_{4,\varepsilon}
h_{\varepsilon}
\nonumber\\
&  \quad -I_{3,\varepsilon}h_{\varepsilon}^2 +  I_{1,\varepsilon}I_{3,\varepsilon}h_{\varepsilon}^2 + I_{3,\varepsilon}I_{4,\varepsilon}g_{\varepsilon}^2+ I_{2,\varepsilon}I_{3,\varepsilon}h_{\varepsilon}^2-2I_{2,\varepsilon}I_{3,\varepsilon}h_{\varepsilon}
g_{\varepsilon}\nonumber\\
&  \quad -  2I_{2,\varepsilon}I_{4,\varepsilon}I_{5,\varepsilon}h_{\varepsilon}g_{\varepsilon}
+2I_{4,\varepsilon}I_{5,\varepsilon}h_{\varepsilon}
g_{\varepsilon}-I_{2,\varepsilon}I_{3,\varepsilon}I_{4,\varepsilon}
g_{\varepsilon}^2+I_{2,\varepsilon}^2I_{3,\varepsilon}h_{\varepsilon}g_{\varepsilon}\\
& \quad +  I_{4,\varepsilon}^2I_{5,\varepsilon}g_{\varepsilon}^2.
\end{align*}
Here
$
g_{\varepsilon}=g_{\varepsilon}\left(\vartheta^*_{\varepsilon},s\right),
$
$
h_{\varepsilon}=h_{\varepsilon}\left(\vartheta^*_{\varepsilon},s\right)
$
and $I_{1,\varepsilon}$, $I_{2,\varepsilon}$, $I_{3,\varepsilon}$, $I_{4,\varepsilon}$, $I_{5,\varepsilon}$ are the ``empirical versions'' of $\hat{I}_1$, $\hat{I}_2$, $\hat{I}_3$, $\hat{I}_4$, $\hat{I}_5$, respectively.

In the construction of the test we introduce one condition else.

\textit{$\mathcal{R}_1.$  We suppose that $\varphi_2\left(r\right),$ $r\in[0,1)$ defined by \eqref{ch112} is strictly positive function.}

We have the uniform convergence in probability w.r.t. $s\in[0,T]$ $(\varepsilon \rightarrow 0)$
$$
\sup_{s\in[0,T]}\left|
\varphi_{2,\varepsilon}(s)- \hat{\varphi}_2(s)\right|\longrightarrow  0.$$
This convergence we obtain due to the consistency of the estimator and the smoothness of the functions $g_\varepsilon\left(\cdot,\cdot\right)$ and $h_\varepsilon\left(\cdot,\cdot\right)$.

Therefore we can introduce the function
$$
\varphi_{2,\varepsilon}^+\left(s\right)=
 \left \{\begin{array}{ll}
  \varphi_{2,\varepsilon}\left(s\right)^{-1}, \quad & \textrm{if}   \;\varphi_{2,\varepsilon}\left(s\right)>0,\\
  0, & \textrm{else},
\end{array}
\right .
$$
which asymptotically coincides with $\hat{\varphi}_2(s)^{-1}$ and therefore the limit distribution does not change.

Hence we consider (formally) the statistic
\begin{align}
\label{ch100}
&W_{\varepsilon}(t) = U_{\varepsilon}\left(\frac{t}{T}\right)\\
& +\frac{1}{T}
\int_0^{t} \int_{0}^s\varphi_{2,\varepsilon}^+(s)[\lambda_1(T)
\varphi_{1,\varepsilon}(s)h_{\varepsilon}(\vartheta^*_{\varepsilon},q)+\lambda_2(T)
  \psi_{2,\varepsilon}(s)
  g_{\varepsilon}(\vartheta^*_{\varepsilon},q)
]\textrm{d}U_{\varepsilon}\left(\frac{q}{T}\right)\textrm{d}s\nonumber
\end{align}
where $\lambda_1(T)=T^2\sqrt{T}, \lambda_2(T)=\displaystyle\frac{1}{T}$ and $U_{\varepsilon}(\cdot)$ was defined by \eqref{ch128}. If we prove that
$$
W_{\varepsilon}(t) \longrightarrow L\left[U\right](\nu)=w_\nu,
$$
then the test based on this statistic will be ADF. The main technical problem in carrying out this program is to define the stochastic integrals
\begin{align}
\label{cch00}
K_{\varepsilon}(\vartheta^*_{\varepsilon},s) = \int_0^s
h_{\varepsilon}\left(\vartheta^*_{\varepsilon},q\right)\;\textrm{d}U_{\varepsilon}\left(\frac{q}{T}\right)
\end{align}
and
\begin{align}
\label{cch01}
L_{\varepsilon}(\vartheta^*_{\varepsilon},s) = \int_0^s
g_{\varepsilon}\left(\vartheta^*_{\varepsilon},q\right)\;\textrm{d}U_{\varepsilon}\left(\frac{q}{T}\right).
\end{align}

Unfortunately we can not calculate them as they are written now, because the integrand contains the MDE $\vartheta^{*}_{\varepsilon}$  and this estimator depends
on the whole trajectory $X^\varepsilon =\left(X_t,0\leq t\leq
T\right)$. Therefore the corresponding stochastic integrals $K_{\varepsilon}\left(\vartheta^*_{\varepsilon},s\right)$ and $L_{\varepsilon}\left(\vartheta^*_{\varepsilon},s\right)$
 are not well defined.

To avoid this problem we use an approach which is based on the application of the It\^{o} formula, i.e., we replace the corresponding stochastic integrals by the ordinary ones. Note that this approach was applied in the similar problem in \cite{Kut13a}.

Introduce the statistic
\begin{align}
\label{ch14}
&K\left(\vartheta,s\right)=\int_0^s\hat{h}\left(\vartheta,q\right)\;\textrm{d}U\left(\frac{q}{T}\right)
=\frac{1}{\sqrt{T}}\int_0^s\hat{h}\left(\vartheta,q\right)\;S\left(\vartheta,x_q\right)\;\textrm{d}u(q),
\end{align}
where the process $U\left(\cdot\right)$ and the function $\hat{h}\left(\cdot,\cdot\right)$ were defined by \eqref{ch1002} and \eqref{cch05}, respectively.
Indeed the It\^{o} formula gives us the following representation
\begin{align*}
&\textrm{d}(\hat{h}(\vartheta,q)S(\vartheta,x_q)u(q))
=\left(\hat{h}'(\vartheta,q)S(\vartheta,x_q)+\hat{h}(\vartheta,q)
S'(\vartheta,x_q)S(\vartheta,x_q)\right)u(q)\textrm{d}q\\
&\qquad\qquad\qquad\qquad\quad\qquad+\hat{h}\left(\vartheta,q\right)S\left(\vartheta,x_q\right)\textrm{d}u(q).
\end{align*}
Here $\hat{h}'\left(\vartheta,q\right)$ is the derivative of $\hat{h}\left(\vartheta,q\right)$ w.r.t. $q$, given by the following expression
\begin{align}
\label{cch04}
\begin{split}
\hat{h}'\left(\vartheta,q\right)=\frac{1}{\sqrt{T}}\;J(\vartheta)^{-1}
\dot{S}'(\vartheta,x_q(\vartheta))\;S\left(\vartheta,x_q(\vartheta)\right).
\end{split}
\end{align}
Therefore the statistic $K\left(\vartheta,\cdot\right)$ defined by \eqref{ch14} can be written as follows:
\begin{align*}
&
K\left(\vartheta,s\right)=\frac{1}{\sqrt{T}}\;\hat{h}\left(\vartheta,s\right)S\left(\vartheta,x_s\right)u(s)\\
&\qquad\qquad-\frac{1}{\sqrt{T}}\int_0^s\left(\hat{h}'\left(\vartheta,q\right)S\left(\vartheta,x_q\right)+\hat{h}\left(\vartheta,q\right)
S'\left(\vartheta,x_q\right)S\left(\vartheta,x_q\right)\right)u(q)\;\textrm{d}q.
\end{align*}
Hence we obtain for the  process \eqref{cch00} the following representation
\begin{align}
\label{cch02}
\begin{split}
&K_{\varepsilon}(\vartheta^*_{\varepsilon},s)=\frac{1}{\sqrt{T}}\;h_\varepsilon\left(\vartheta^*_{\varepsilon},s\right)
S\left(\vartheta^*_{\varepsilon},X_s\right)u_\varepsilon(s)\\
&-\frac{1}{\sqrt{T}}\int_0^s\left(h_\varepsilon'(\vartheta^*_{\varepsilon},q)
S(\vartheta^*_{\varepsilon},X_q)+h_\varepsilon(\vartheta^*_{\varepsilon},q)
S'(\vartheta^*_{\varepsilon},X_q)S(\vartheta^*_{\varepsilon},X_q)
\right)u_\varepsilon(q)\;\textrm{d}q
\end{split}
\end{align}
and the integral is now well defined. Here
$$
h'_\varepsilon\left(\vartheta^*_{\varepsilon},q\right)=\frac{1}{\sqrt{T}}\;J_\varepsilon(\vartheta^*_{\varepsilon})^{-1}
\dot{S}'(\vartheta^*_{\varepsilon},X_q)\;S\left(\vartheta^*_{\varepsilon},X_q\right).
$$

Similarly,  the process \eqref{cch01} can be written as follows:
\begin{align}
\label{ch03}
\begin{split}
&L_{\varepsilon}(\vartheta^*_{\varepsilon},s)=\frac{1}{\sqrt{T}}\;g_\varepsilon\left(\vartheta^*_{\varepsilon},s\right)
S\left(\vartheta^*_{\varepsilon},X_s\right)u_\varepsilon(s)\\
&-\frac{1}{\sqrt{T}}\int_0^s\left(g_\varepsilon'(\vartheta^*_{\varepsilon},q)
S(\vartheta^*_{\varepsilon},X_q)+g_\varepsilon(\vartheta^*_{\varepsilon},q)
S'(\vartheta^*_{\varepsilon},X_q)S(\vartheta^*_{\varepsilon},X_q)
\right)u_\varepsilon(q)\textrm{d}q,
\end{split}
\end{align}
where
$$
g'_\varepsilon(\vartheta^*_{\varepsilon},q)=-S(\vartheta^*_{\varepsilon},X_q)^{-1}
S'(\vartheta^*_{\varepsilon},X_q)\int_{q}^TS(\vartheta^*_{\varepsilon},X_s)\;
\dot{x}_s(\vartheta^*_{\varepsilon})\;\textrm{d}s-\dot{x}_q(\vartheta^*_{\varepsilon})
$$
is the ``empirical version'' of
\begin{align}
\begin{split}
\label{cch03}
&\hat{g}'(\vartheta,q)=-S(\vartheta,x_q(\vartheta))^{-1}
S'(\vartheta,x_q(\vartheta))\int_{q}^TS(\vartheta,x_s(\vartheta))\;\dot{x}_s(\vartheta)\;\textrm{d}s-\dot{x}_q(\vartheta).
\end{split}
\end{align}

Then the formal expression \eqref{ch100} for $W_\varepsilon(t)$ can be replaced by
\begin{align}
\label{cch09}
\begin{split}
&\tilde{W}_{\varepsilon}(t) = U_{\varepsilon}\left(\frac{t}{T}\right)\\
& +\frac{1}{T}
\displaystyle\int_0^{t} \varphi_{2,\varepsilon}^+\left(s\right)\left[\lambda_1(T)\;
\varphi_{1,\varepsilon}(s)\;K_{\varepsilon}(\vartheta^*_{\varepsilon},s)+\lambda_2(T)\;
  \psi_{2,\varepsilon}(s)\;
  L_{\varepsilon}(\vartheta^*_{\varepsilon},s)
\right]\;\textrm{d}s,
\end{split}
\end{align}
where $\lambda_1(T)=T^2\sqrt{T}, \lambda_2(T)=\displaystyle\frac{1}{T}$ and the processes $U_{\varepsilon}(\cdot)$, $K_{\varepsilon}\left(\vartheta^*_{\varepsilon},\cdot\right)$ and $L_{\varepsilon}\left(\vartheta^*_{\varepsilon},\cdot\right)$ admit the representations \eqref{ch128}, \eqref{cch02} and \eqref{ch03}, respectively.

The test is given in the following Theorem.
\begin{theorem}
\label{ch1T2}
Suppose that the conditions of regularity $\mathcal{R}$  and $\mathcal{R}_1$ are fulfilled, then the test
$$
\psi_{\varepsilon} = \1_{\left\{ \Delta_\varepsilon> c_\alpha \right\}},
\qquad\Pb\left(\tilde{\Delta}> c_{\alpha}\right)=\alpha,
$$
with
\begin{align*}
\Delta_{\varepsilon} = \frac{1}{T}\int_0^T
\tilde{W}_{\varepsilon}(t)^2  \textrm{d}t,\quad\tilde{\Delta} \equiv
\int_0^1 w_{\nu}^2\;\textrm{d}\nu
\end{align*}
is ADF and of asymptotic size $\alpha\in(0,1)$.
\end{theorem}
\begin{pf}
We have to show that, under hypothesis $\mathcal{H}_0$, the convergence
\begin{equation}
\label{ch131}
\Delta_{\varepsilon} \Longrightarrow \tilde{\Delta}
\end{equation}
holds.

Recall that the process $U_{\varepsilon}(\cdot)$ has the following representation
\begin{align*}
&U_{\varepsilon}\left(\frac{t}{T}\right) =  \frac{1}{\sqrt{T}}\;S\left(\vartheta^{*}_{\varepsilon},X_t\right)\;u_\varepsilon\left(t\right)\\
&\qquad \qquad\qquad- \frac{1}{\sqrt{T}}
\int_0^tS'\left(\vartheta^{*}_{\varepsilon},X_s\right)\;S\left(\vartheta^{*}_{\varepsilon},X_s\right)
\;u_\varepsilon\left(s\right)\;\textrm{d}s.
\end{align*}
Note that we have already the convergence in probability uniformly w.r.t. $s\in[0,T]$ (as $\varepsilon \rightarrow 0$)
$$
\sup_{s\in[0,T]}\left|X_s
- x_s(\vartheta)\right|\longrightarrow 0, \qquad\sup_{s\in[0,T]}\left|u_\varepsilon\left(s\right)
- u\left(s\right)\right|\longrightarrow 0.
$$
Further, we can write
\begin{align*}
&\left|S\left(\vartheta^{*}_{\varepsilon},X_s\right)-S\left(\vartheta,x_s\right)\right|
\leq \left|S\left(\vartheta^{*}_{\varepsilon},X_s\right)-S\left(\vartheta,X_s\right)\right|+ \left|S\left(\vartheta,X_s\right)-S\left(\vartheta,x_s\right)\right|\\
&\qquad\qquad\qquad\qquad\qquad\leq \left|\vartheta^{*}_{\varepsilon}-\vartheta\right|\left|\dot{S}(\tilde{\vartheta},X_s)\right|
+\left|X_s-x_s\right|\left|S'(\vartheta,\tilde{X}_s)\right|.
\end{align*}
Here
$|\tilde{\vartheta} - \vartheta| \leq |\vartheta^{*}_{\varepsilon}- \vartheta|$
and
$$
\left|\tilde{X}_s - X_s\right|\leq \left| x_s(\vartheta^*_{\varepsilon}) - X_s\right| \leq \left| x_s(\vartheta^*_{\varepsilon})-x_s(\vartheta)\right|+\left|x_s(\vartheta)-X_s\right|\rightarrow 0.
$$
This convergence is uniform w.r.t. $s\in[0,T]$.

Then we know that the functions $\dot{S}(\vartheta,x)$ and $S'(\vartheta,x)$ are bounded by regularity conditions $\mathcal{R}$, the process $X_s$ converges uniformly w.r.t. $s\in[0,T]$ to $x_s\left(\vartheta\right)$ and due to the consistency of the estimator $\vartheta^{*}_{\varepsilon}$, we obtain the uniform convergence w.r.t. $s\in[0,T]$ (in probability)
$$
\sup_{s\in[0,T]}\left|S\left(\vartheta^{*}_{\varepsilon},X_s\right)-S\left(\vartheta,x_s\right)\right|\longrightarrow 0.
$$
Further, similar arguments give the uniform convergence w.r.t. $s\in[0,T]$ (in probability)
$$
\sup_{s\in[0,T]}\left|S'\left(\vartheta^{*}_{\varepsilon},X_s\right)-S'\left(\vartheta,x_s\right)\right|\longrightarrow 0.
$$
Therefore, we obtain the uniform convergence w.r.t. $t\in[0,T]$ (in probability)
\begin{align*}
&U_{\varepsilon}\left(\frac{t}{T}\right) \longrightarrow  \frac{1}{\sqrt{T}}\;S(\vartheta,x_t)\;u(t)- \frac{1}{\sqrt{T}}
\int_0^tS'(\vartheta,x_s)\;S(\vartheta,x_s)\;u(s)\;\textrm{d}s\\
&\qquad \qquad \quad = \frac{1}{\sqrt{T}}\int_0^tS(\vartheta,x_s)\;\textrm{d}u(s)= U\left(\frac{t}{T}\right).
\end{align*}

Now we have to show that
$K_{\varepsilon}\left(\vartheta^*_{\varepsilon},s\right)\longrightarrow K(\vartheta,s)$, where
\begin{align*}
&K_{\varepsilon}(\vartheta^*_{\varepsilon},s)=\frac{1}{\sqrt{T}}\;h_\varepsilon\left(\vartheta^*_{\varepsilon},s\right)
S\left(\vartheta^*_{\varepsilon},X_s\right)u_\varepsilon(s)\\
&-\frac{1}{\sqrt{T}}\int_0^s\left(h_\varepsilon'(\vartheta^*_{\varepsilon},q)
S(\vartheta^*_{\varepsilon},X_q)+h_\varepsilon(\vartheta^*_{\varepsilon},q)
S'(\vartheta^*_{\varepsilon},X_q)S(\vartheta^*_{\varepsilon},X_q)
\right)u_\varepsilon(q)\;\textrm{d}q
\end{align*}
and
\begin{align*}
&
K\left(\vartheta,s\right)=\frac{1}{\sqrt{T}}\;\hat{h}\left(\vartheta,s\right)S\left(\vartheta,x_s\right)u(s)\\
&\qquad\quad\quad-\frac{1}{\sqrt{T}}\int_0^s\left(\hat{h}'\left(\vartheta,q\right)S\left(\vartheta,x_q\right)+\hat{h}\left(\vartheta,q\right)
S'\left(\vartheta,x_q\right)S\left(\vartheta,x_q\right)\right)u(q)\;\textrm{d}q,
\end{align*}
where $\hat{h}\left(\cdot,\cdot\right)$ and $\hat{h}'(\cdot,\cdot)$ were defined by \eqref{cch05}$-$\eqref{cch04}.

We have
\begin{align*}
&\left|\dot{S}\left(\vartheta^*_{\varepsilon},X_s\right)-
\dot{S}\left(\vartheta,x_s\right)\right|\leq \left|\dot{S}\left(\vartheta^*_{\varepsilon},X_s\right)-
\dot{S}\left(\vartheta,X_s\right)\right| +\left|\dot{S}\left(\vartheta,X_s\right)-\dot{S}\left(\vartheta,x_s\right)\right|\\
&\qquad \qquad \qquad\qquad\quad\quad\;\leq\left|\vartheta^*_{\varepsilon}-\vartheta\right|\left|\ddot{S}(\tilde{\vartheta},X_s)\right|
+\left|X_s-x_s\right|\left|\dot{S}'(\vartheta,\tilde{X}_s)\right|\\
&\qquad \qquad \qquad\qquad\quad\quad\;\leq \tilde{C}_1\left|\vartheta^*_{\varepsilon}-\vartheta\right|
+\tilde{C}_2\left|X_s-x_s\right|.
\end{align*}
Here $\tilde{C}_1$ and $\tilde{C}_2$ are some constants and $\ddot{S}(\vartheta,x)$, which means the second derivative w.r.t. $\vartheta$, and $\dot{S}'(\vartheta,x)$ are bounded functions due to the regularity conditions $\mathcal{R}$. Therefore we have the convergence in probability
$$
\sup_{s\in[0,T]}\left|\dot{S}\left(\vartheta^{*}_{\varepsilon},X_s\right)-\dot{S}\left(\vartheta,x_s\right)\right|\longrightarrow 0.
$$

 Further, we have
$$
J_\varepsilon(\vartheta^*_{\varepsilon})- J(\vartheta)=
\displaystyle\int_0^T\dot{x}_s(\vartheta^*_{\varepsilon})^2\textrm{d}s -
\displaystyle\int_0^T\dot{x}_s(\vartheta)^2\textrm{d}s\longrightarrow 0.
$$
Hence we obtain the convergence in probability
$$
\sup_{s\in[0,T]}\left|h_{\varepsilon}\left(\vartheta^*_{\varepsilon},s\right)-
\hat{h}\left(\vartheta,s\right)\right|\longrightarrow
0.
$$
Similarly, it is shown that we have the convergence in probability
$$
\sup_{s\in[0,T]}\left|h'_{\varepsilon}\left(\vartheta^*_{\varepsilon},s\right)-
\hat{h}'\left(\vartheta,s\right)\right|\longrightarrow
0.
$$
Therefore we obtain the convergence in probability, by the uniform convergence of $u_\varepsilon\left(s\right)$ to $u(s)$ w.r.t. $s\in[0,T]$,
$$
K_{\varepsilon}\left(\vartheta^*_{\varepsilon},s\right)\longrightarrow K(\vartheta,s).
$$

Now we have to show that
$L_{\varepsilon}(\vartheta^*_{\varepsilon},s) \longrightarrow L(\vartheta,s),$ where
\begin{align*}
&L_{\varepsilon}(\vartheta^*_{\varepsilon},s)=\frac{1}{\sqrt{T}}\;g_\varepsilon\left(\vartheta^*_{\varepsilon},s\right)
S\left(\vartheta^*_{\varepsilon},X_s\right)u_\varepsilon(s)\\
&-\frac{1}{\sqrt{T}}\int_0^s\left(g_\varepsilon'(\vartheta^*_{\varepsilon},q)
S(\vartheta^*_{\varepsilon},X_q)+g_\varepsilon(\vartheta^*_{\varepsilon},q)
S'(\vartheta^*_{\varepsilon},X_q)S(\vartheta^*_{\varepsilon},X_q)
\right)u_\varepsilon(q)\;\textrm{d}q
\end{align*}
and
\begin{align*}
&
L\left(\vartheta,s\right)=\frac{1}{\sqrt{T}}\;\hat{g}\left(\vartheta,s\right)S\left(\vartheta,x_s\right)u(s)\\
&\qquad\qquad-\frac{1}{\sqrt{T}}\int_0^s\left(\hat{g}'\left(\vartheta,q\right)S\left(\vartheta,x_q\right)+\hat{g}\left(\vartheta,q\right)
S'\left(\vartheta,x_q\right)S\left(\vartheta,x_q\right)\right)u(q)\;\textrm{d}q,
\end{align*}
where $\hat{g}\left(\cdot,\cdot\right)$ and $\hat{g}'(\cdot,\cdot)$ were defined by \eqref{cch06} and \eqref{cch03}, respectively.

Observe that for $s\in[0,T]$, we have
\begin{align*}
&\left|I_{\varepsilon}(\vartheta^*_{\varepsilon},s)
- I(\vartheta,s)\right|
\leq
\int_s^T \left|S(\vartheta^*_{\varepsilon},X_v)
\left(\dot{x}_v(\vartheta^*_{\varepsilon})-\dot{x}_v(\vartheta)\right)\right|\textrm{d}v \\ &\qquad\qquad\qquad\qquad\qquad+\int_s^T \left|\dot{x}_v(\vartheta)
\left(S(\vartheta^*_{\varepsilon},X_v)
-S(\vartheta,x_v(\vartheta))\right)\right|\textrm{d}v\\
&\qquad\qquad\qquad\qquad\qquad\leq
\int_s^T \left|S(\vartheta^*_{\varepsilon},X_v)\right|
\left|\dot{x}_v(\vartheta^*_{\varepsilon})-\dot{x}_v(\vartheta)\right|\textrm{d}v \\ &\qquad\qquad\qquad\qquad\qquad+\left|\vartheta^*_{\varepsilon}-\vartheta\right| \int_s^T \left|\dot{x}_v(\vartheta)\right| \left|\dot{S}(\tilde{\vartheta},X_v)\right|\textrm{d}v\\
&\qquad\qquad\qquad\qquad\qquad+\int_s^T \left|\dot{x}_v(\vartheta)\right| \left|X_v-x_v(\vartheta)\right| \left|S'(\vartheta,\tilde{X}_v)\right|\textrm{d}v\\
&\qquad\qquad\leq \left(\tilde{C}_3\left|\vartheta^*_{\varepsilon}-\vartheta\right|
+\tilde{C}_4\left|\vartheta^*_{\varepsilon}-\vartheta\right|+\tilde{C}_5\sup_{s\in[0,T]}\left|X_s-x_s(\vartheta)\right|\right)\left(T-s\right),
\end{align*}
where $\tilde{C}_3$, $\tilde{C}_4$ and $\tilde{C}_5$ are constants. Therefore we obtain the  convergence in probability uniformly w.r.t. $s\in[0,T]$
$$
\sup_{s\in[0,T]}\left|I_{\varepsilon}(\vartheta^*_{\varepsilon},s)
- I(\vartheta,s)\right|\longrightarrow 0$$
because the estimator $\vartheta^*_{\varepsilon}$ is consistent, the process $X_s$ converges uniformly w.r.t. $s\in[0,T]$ to $x_s(\vartheta)$ and the derivatives are bounded due to the conditions of regularity $\mathcal{R}$. Further, we proved already the uniform convergence w.r.t. $s\in[0,T]$ of $u_\varepsilon(s)$ to $u(s)$ and $S(\vartheta^*_{\varepsilon},X_s)$ to $S(\vartheta,x_v)$. Hence, we have the convergence in probability
$$
\sup_{s\in[0,T]}\left|g_{\varepsilon}\left(\vartheta^*_{\varepsilon},s\right)-
\hat{g}\left(\vartheta,s\right)\right|\longrightarrow 0.$$
Similarly, by the regularity conditions $\mathcal{R}$, it is shown that (in probability)
$$
\sup_{s\in[0,T]}\left|g'_{\varepsilon}\left(\vartheta^*_{\varepsilon},s\right)-
\hat{g}'\left(\vartheta,s\right)\right|\longrightarrow 0.
$$

Therefore we obtain the convergence in probability
$$
L_{\varepsilon}\left(\vartheta^*_{\varepsilon},s\right)
\longrightarrow L\left(\vartheta,s\right).
$$
Further, a similar arguments give the convergence  in probability uniformly w.r.t. $s\in[0,T]$
due to the regularity conditions $\mathcal{R}$ and the consistency of the estimator $\vartheta^*_{\varepsilon}$
$$\sup_{s\in[0,T]}\left| \psi_{2,\varepsilon}(s) -
\hat{\psi}_2(s)\right|\longrightarrow  0
$$
and
$$
\sup_{s\in[0,T]}\left| \varphi_{1,\varepsilon}(s) -
\hat{\varphi}_1(s)\right|\longrightarrow  0, \sup_{s\in[0,T]}\left|
\varphi_{2,\varepsilon}(s)- \hat{\varphi}_2(s)\right|\longrightarrow  0.$$

Finally, the convergence mentioned in
\eqref{ch131} is proved and using this result the test
$\psi_{\varepsilon}$ is ADF and of asymptotic size $\alpha \in (0,1)$.
\end{pf}

\section{The case of MLE}
\label{ch1s05}

This case was studied in \cite{Kh81}$-$\cite{KK14}. They proposed a
linear transformation, which yields the convergence of the test statistic to
the integral of Wiener process. Therefore they showed that the test based on
this statistic is ADF. To obtain the linear transformation mentioned in
\cite{KK14}, we put $h(\vartheta,r)=g(\vartheta,r)$ in \eqref{ch110} and
we obtain
$$
I_1 = I_2 =I_4= \int_0^r h(\vartheta,q)^2\textrm{d}q,\quad I_3 =I_5=\int_0^r h(\vartheta,q)\;\textrm{d}q.
$$
 Therefore we can write
\begin{align*}
\varphi_1(r) & =  h-h+I_1h-3I_1h+I_3h^2+  I_3h^2+2I_1h-I_1^2h-2I_1I_3h^2-I_1^2I_3h^2\\
& \quad - I_1h+I_1^2I_3h^2 - I_1^2h+3I_1^2h+  I_1I_3h^2-2I_1I_3h^2-2I_1^2h+  2I_1I_3h^2\\
&  \quad +  I_1^2I_3h^2+I_1^3h+I_3h^2+I_1h+ 2I_1^2h-I_1^2h + I_1^3h + I_1^2I_3h^2-  2I_1^2I_3h^2\\
&  \quad -I_1^3h+I_1I_3h^2+I_1I_3h^2+2I_1I_3h^2-I_1^3h
-  I_1^2I_3h^2-I_3h^2-2I_1I_3h^2\\
 & \quad  +I_1^2I_3h^2-I_1I_3h^2
- 2I_3h^2\\
& =  0,\\
\varphi_2(r) & = 1 + I_3h-3I_1I_3h+I_1I_3h+I_3h-3I_1I_3h+I_1I_3h-I_3h+6I_1^2-2I_1^4\\
&  \quad+   3I_1^2I_3h- 2I_1^2I_3h+2I_1I_3h-I_1^3I_3h+I_1^3I_3h-I_1^3I_3h +  I_1^2-4I_1 \\
&  \quad +  I_1I_3h-I_1^2I_3h+2I_1^2I_3h+I_1I_3h-I_1^2I_3h +  2I_1^3 +  I_1^2I_3h-  4I_1^3\\
&  \quad+  I_1^3I_3h-I_1^3I_3h-I_1^2I_3h + I_1^3I_3h- I_1^2I_3h+4I_1^3+ 3I_1^2I_3h-  I_1^3I_3h\\
&  \quad+   I_1^4-2I_1^3+I_1^4+2I_1- 2I_1^2+I_1^3I_3h - 2I_1^2I_3h-I_1I_3h - I_1^2I_3h-4I_1^2 \\
& =  1+I_3h-I_1I_3h+I_1^2-2I_1\\
& = (1-I_1)(1+I_3h-I_1)
\end{align*}
and
\begin{align*}
\psi_2(r) & = h+I_3h^2+I_1h+I_3h^2-3I_1h-2I_1I_3h^2-2I_1^2h-2I_1I_3h^2-I_1^3h  \\
&  \quad +  I_1^3h+I_1^2I_3h^2-I_1^3h-I_1I_3h^2+I_1^2h+I_1h -  I_1^2h-I_3h^2+  3I_1^2h \\
&  \quad -  I_1^2h+I_1^3h+I_1I_3h^2 + I_1I_3h^2+ I_1I_3h^2-I_1^2I_3h^2+2 I_1^2I_3h^2\\
&  \quad -  I_1^2I_3h^2-2I_1^2I_3h^2 +2I_1I_3h^2  +I_1^2I_3h^2\\
& = h\;(1+I_3h-I_1),
\end{align*}
where $h=h(\vartheta,r)$. Hence the linear transformation \eqref{ch114} will have the following expression:
\begin{equation}
\label{ch132}
L[U](\nu) = U(\nu)+\displaystyle\int_0^{\nu} \int_0^r h(r)\;\mathbb{N}(r)^{-1}h(q)\;\textrm{d}U(q)\;\textrm{d}r=w_{\nu},
\end{equation}
where
 $\mathbb{N}(r) = \displaystyle\int_r^1 h(q)^2\textrm{d}q$
 and
 $$
 U(\nu) = W(\nu) -
\int_0^1h(r)\;\textrm{d}W(r)\;\int_0^\nu h(r)\;\textrm{d}r,
$$
 with $W(\nu)$ and $w_{\nu}, 0\leq\nu \leq 1$ are some standard Wiener processes.

The  transformation $L[U](\cdot)$ of the limit process $U(\cdot)$
given by \eqref{ch132} coincides with one by Khmaladze \cite{Kh81}.

{\bf Acknowledgements.}
 The author is grateful to Yu. A. Kutoyants for the statement of the problem and many helpful suggestions during the preparation of this paper. The author would like to thank the two Referees for the comments which allowed to improve essentially the exposition.



\end{document}